\documentclass[12pt]{amsart}

\usepackage{hyperref}
\usepackage{amssymb, latexsym, amsmath, amsfonts}
\usepackage{amscd}
\usepackage{curves}
\usepackage{amsthm}
\usepackage{amsbsy}

\newcommand{\myamstitle}[3]{
  \title[#2]{#1}
  \author{Alexandru Scorpan}
  \address{Department of Mathematics, University of California, Berkeley\\ 970 Evans Hall, Berkeley, CA 94720}
  \email{scorpan@math.berkeley.edu}
  \urladdr{www.math.berkeley.edu/\textasciitilde scorpan}
  \date{#3}
  }

\theoremstyle{plain} 

\newtheorem{theorem}{Theorem}[section] 
\newtheorem{lemma}[theorem]{Lemma}
\newtheorem{proposition}[theorem]{Proposition}
\newtheorem{corollary}[theorem]{Corollary}

\newtheorem{question}[theorem]{Question}

\newtheorem{conjecture}[theorem]{Conjecture}

\theoremstyle{definition} 

\theoremstyle{remark} 

\newtheorem{remark}[theorem]{Remark}

\newtheorem{example}[theorem]{Example}
\newtheorem{piece}[theorem]{}

   {\endgroup\medskip}

\newcommand{\proofbox}{\ensuremath{\smash{\Box}}}
\newcommand{\noproof}{\mbox{}\hfill\proofbox}     

\newenvironment{Proof}[1]{\begin{proof}[Proof of #1]}{\end{proof}} 

\newcommand{\Secref}[1]{Sec\-tion \ref{#1}}
\newcommand{\Lemmaref}[1]{Lem\-ma \ref{#1}}
\newcommand{\Propref}[1]{Pro\-po\-si\-tion \ref{#1}} 
\newcommand{\Thmref}[1]{Theo\-rem \ref{#1}}

\newcommand{\Corref}[1]{Co\-rol\-la\-ry \ref{#1}}
\newcommand{\Conjref}[1]{Con\-jec\-ture \ref{#1}}
\newcommand{\Figref}[1]{Fi\-gu\-re \ref{#1}} 
\newcommand{\Exref}[1]{Ex\-am\-ple \ref{#1}}

\newcommand{\mraise}[2]{\raisebox{#1}{$#2$}}      

\newcommand{\sub}[1]{\raisebox{-2pt}{$\!_{#1}$} }  

\newcommand{\adjustnabla}[1]{\sub{#1}} 

\newcommand{\defemph}[1]{{\sffamily\slshape #1}}  

\renewcommand{\phi}{\varphi}
\renewcommand{\theta}{\vartheta}
\renewcommand{\epsilon}{\varepsilon} 

\newcommand{\ie}{{\it i.e.}~}
\newcommand{\eg}{{\it e.g.}~}
\newcommand{\st}{such that}
\newcommand{\Iff}{if and only if}

\newcommand{\LC}{Levi-Civit\`a}
\newcommand{\SW}{Sei\-berg--Witten}

\newcommand{\ac}{al\-most-com\-plex}
\newcommand{\riem}{Riemannian}

\newcommand{\str}{struc\-ture}

\newcommand{\aR}{\mathbb{R}}
\newcommand{\Zi}{\mathbb{Z}}

\newcommand{\Ce}{\mathbb{C}}

\newcommand{\Proj}{\mathbb{P}}

\newcommand{\CP}{\Ce\Proj}

\newcommand{\Sph}[1]{{\mathbb S}^{#1}}  

\newcommand{\del}{\partial}           
\renewcommand{\Bar}[1]{\overline{#1}} 
\newcommand{\rec}[1]{\tfrac{1}{#1}}   
\newcommand{\iso}{\approx}            

\newcommand{\maps}{\longmapsto}       
\newcommand{\longto}{\longrightarrow} 

\newcommand{\inner}[1]{\langle #1\rangle} 
\newcommand{\Inner}[1]{\bigl\langle #1\bigr\rangle} 

\newcommand{\module}[1]{\left|#1\right|}  

\newcommand{\rest}[1]{|_{#1}}                          

\newcommand{\T}[1]{T_{#1}}             
\newcommand{\N}[1]{N_{#1}}             

\newcommand{\Nabla}[1]{\nabla\adjustnabla{#1}}

\newcommand{\pinc}{pin$^{\!{\bf C}}$} 
\newcommand{\spinc}{s\pinc}         

\newcommand{\tsum}{\operatorname{\mraise{1pt}{\sum}}} 

\newcommand{\F}{{\mathcal F}}
\newcommand{\G}{{\mathcal G}}

\newcommand{\comment}[1]{%
	#1
	}

\begin{document}

\myamstitle{A quick survey of foliations\\ on $4$-manifolds}{A quick survey of foliations on $4$-manifolds}{February 7, 2003 (revised \today)}


\begin{abstract}
We present a few general results on foliations of $4$-manifolds by surfaces: existence, tautness, relations to minimal genus of embedded surfaces;
as well as some open problems.
We hope to stimulate interest in this area.
\end{abstract}

\maketitle

\footnotesize

\tableofcontents

\normalsize

\section{Introduction}

While there is a wealth of known facts concerning codimension-$1$ foliations (and in particular foliations on $3$-manifolds), the geometric results are quite scarce once one considers higher codimensions, for example $2$-dimensional foliations on $4$-manifolds. And while there are strong suggestions that the field might yield remarkable insights, one often accuses a lack of methods for attacking it. 

This survey has two goals: to quickly gather the known (geometric-topologic) results concerning $2$-dimensional foliations on $4$-manifolds, and also to frame the context of possible future developments. The biases of the author will become obvious soon enough.

\subsection{Quick definition}

\begin{figure}[bth]
\setlength{\unitlength}{1.2pt}
\begin{picture}(200,130)(-15,0)

\put(0,0){\line(1,0){140}}
\put(0,5){\line(1,0){140}}
\put(0,10){\line(1,0){140}}
\put(0,15){\line(1,0){140}}
\put(0,20){\line(1,0){140}}
\put(140,0){\line(1,2){30}}
\put(140,5){\line(1,2){30}}
\put(140,10){\line(1,2){30}}
\put(140,15){\line(1,2){30}}
\put(140,20){\line(1,2){30}}

\put(0,45){\line(1,0){140}}
\put(0,25){\line(1,0){140}}
\put(0,30){\line(1,0){140}}
\put(0,35){\line(1,0){140}}
\put(0,40){\line(1,0){140}}
\put(140,45){\line(1,2){30}}
\put(140,25){\line(1,2){30}}
\put(140,30){\line(1,2){30}}
\put(140,35){\line(1,2){30}}
\put(140,40){\line(1,2){30}}

\put(0,50){\line(1,0){140}}
\put(0,55){\line(1,0){140}}
\put(0,60){\line(1,0){140}}
\put(0,65){\line(1,0){140}}
\put(0,70){\line(1,0){140}}
\put(140,50){\line(1,2){30}}
\put(140,55){\line(1,2){30}}
\put(140,60){\line(1,2){30}}
\put(140,65){\line(1,2){30}}
\put(140,70){\line(1,2){30}}

\put(0,70){\line(1,2){30}}
\put(30,130){\line(1,0){140}}

\end{picture}
\caption{Local model for a foliation}
\label{fig-foliation.loc}
\end{figure}
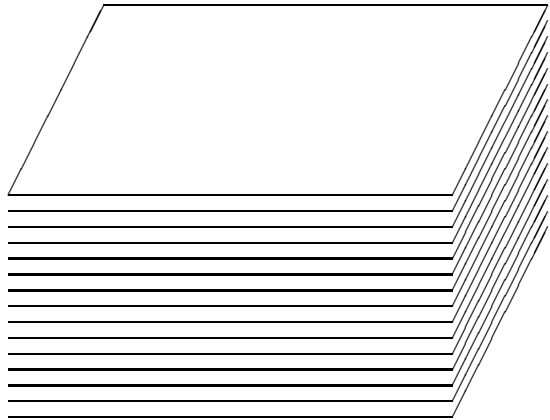

In general, a $k$-dimensional \defemph{foliation} $\F$ is a partition of a manifold $M^m$ into ``slices'' called \defemph{leaves}.
Locally, a foliation must look like the partition of $\aR^m$ into $\aR^k\times\{0\}$ (see \Figref{fig-foliation.loc}).

\begin{example}
An important example is the \defemph{Reeb foliation} of the solid torus. The boundary torus is a leaf, and the interior is foliated by leaves that are planes. (See \Figref{fig-reeb}.)
\end{example}

\comment{
\begin{figure}[bth]
\setlength{\unitlength}{.5pt}
\begin{picture}(500,415)(-250,-265)


\linethickness{1pt}

\qbezier(100,0)(150,0)(150,-100)
\qbezier(150,-100)(150,-200)(100,-200)
\qbezier(100,-200)(83,-200)(75,-195)
\qbezier(100,0)(75,0)(65,-29)
\qbezier(-100,0)(-150,0)(-150,-100)
\qbezier(-150,-100)(-150,-200)(-100,-200)
\qbezier(-100,-200)(-75,-200)(-65,-177)
\qbezier(-100,0)(-67,0)(-60,-55)

\qbezier(100,-200)(260,-200)(250,0)
\qbezier(250,0)(240,150)(0,150)
\qbezier(-100,-200)(-260,-200)(-250,0)
\qbezier(-250,0)(-240,150)(0,150)
\qbezier(100,0)(130,0)(125,35)
\qbezier(125,35)(120,70)(0,70)
\qbezier(-100,0)(-130,0)(-125,35)
\qbezier(-125,35)(-120,70)(0,70)


\linethickness{.5pt}

\qbezier(60,-125)(60,-25)(20,-25)
\qbezier(60,-125)(60,-200)(20,-200)
\qbezier(20,-25)(-5,-25)(-10,-60)
\qbezier(20,-200)(0,-200)(-7,-180)

\qbezier(33,-28)(60,-45)(115,-17)
\qbezier(24,-200)(75,-200)(125,-165)
\qbezier(-120,-40)(-100,-60)(-25,-65)
\qbezier(-25,-175)(-75,-175)(-120,-160)
\qbezier(-25,-65)(40,-75)(50,-105)
\qbezier(50,-105)(65,-175)(-25,-175)


\thinlines
\qbezier(-180,-250)(-220,-200)(-180,-150)
\thicklines
\put(-180,-150){\vector(2,3){0}}
\put(-180,-260){\makebox(0,0){\tt torus leaf}}

\thinlines
\qbezier(80,-250)(120,-200)(80,-150)
\thicklines
\put(80,-150){\vector(-2,3){0}}
\put(80,-260){\makebox(0,0){\tt plane leaf}}

\end{picture}
\caption{Reeb foliation}
\label{fig-reeb}
\end{figure}
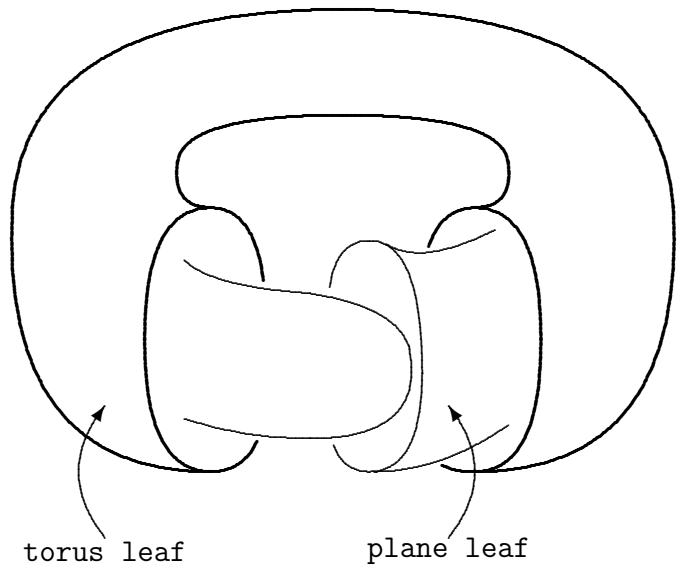
}

\noindent
For a rigorous introduction to foliations, see for example \cite{candelconlon}.
For a general survey of $4$-manifold theory, see \cite{gompf}.

\medskip

Unless otherwise specified, all  foliations considered from now on are oriented $2$-dimensional foliations of oriented closed $4$-manifolds.

\subsection{Why bother?}

One hint that foliations on $4$-manifolds are worth studying (especially \emph{taut} foliations) is Kronheimer's \Thmref{thm-kron} and \Conjref{conj}. Taut foliations might offer \emph{minimal genus bounds} for embedded surfaces.

In a slightly larger context, the relationship (if any) between foliations on $4$-manifolds and \emph{\SW\ theory} is worth elucidating.

Another question worth asking is: For what foliations is the induced \ac\ \str\ ``nice'' (\ie close to symplectic). One such problem asks for which foliations does the induced \ac\ \str\ have \emph{Gromov compactness} (\ie whether the space of $J$-holomorphic curves of a fixed genus and homology class is compact). 

In general, one can hope that foliations will help better visualize, manipulate and understand \emph{\ac\ \str s}, maybe in a manner similar to the one in which open-book decompositions help understand contact \str s on $3$-manifolds.

\subsection{Larger context}

\begin{figure}[bth]
\setlength{\unitlength}{0.9pt}
\begin{picture}(390,335)(-195,-160)

\newcommand{\C}[1]{\makebox(0,0){#1}}
\newcommand{\B}[1]{\oval(100,20)\C{#1}}
\newcommand{\f}{\footnotesize}
\newcommand{\x}{\scriptsize\bf}
\thicklines

\thinlines
\put(-195,-160){\line(1,0){390}}
\put(-195,-160){\line(0,1){335}}
\put(195,175){\line(-1,0){390}}
\put(195,175){\line(0,-1){335}}
\thicklines

\put(-70,40){\B{\spinc-structure}}
\put(70,40){\B{almost-complex}}
\put(0,-30){\B{foliation}}

\put(-120,140){\B{Seiberg--Witten}}
\put(120,140){\B{symplectic}}
\put(0,-140){\B{taut foliation}}

\put(-140,-80){\B{genus bound}}
\put(140,-80){\B{Lefschetz pencil}}

\thicklines
\put(20,-20){\vector(1,2){24.5}} 
\put(20,40){\vector(-1,0){39}} 
\put(120,130){\vector(-1,-2){39.3}} 
\put(0,-130){\vector(0,1){89}} 
\put(-80,50){\vector(-1,2){39.5}} 

\qbezier(145,-90)(110,-135)(54.5,-140.2) 
\put(50.5,-140){\vector(-1,0){0}}

\qbezier(-94,151.8)(0,185)(95,150.5)
\put(-97,151){\vector(-3,-1){0}}
\put(-47,154){\C{\f Taubes}}

\qbezier(-145,130)(-220,40)(-161.2,-66)
\put(-160,-69){\vector(1,-2){0}}
\put(-111,-55){\C{\f Kronheimer--Mrowka}}

\qbezier(144,125.5)(220,40)(162,-66)
\put(160,-69.5){\vector(-1,-2){0}}
\put(140,129.5){\vector(-1,1){0}}
\put(137,-55){\C{\f Donaldson}}
\put(170,120){\C{\f Gompf}}

\thinlines
\qbezier(-143,-94)(-110,-135)(-50,-140)
\thicklines
\put(-145,-90.5){\vector(-2,3){0}}
\put(-80,-130){?}
\put(-103,-100){\C{\f Kronheimer}}

\qbezier(65,30)(70,0)(50,-20) 
\put(47.5,-22.5){\vector(-1,-1){0}} 

\end{picture}
\caption{Context}
\label{fig-context}
\end{figure}

In \Figref{fig-context} we have tried to map the larger mathematical neighborhood in which $2$-dimensional foliations on $4$-manifolds could be placed:

\defemph{Foliations}, in the presence of a metric, induce \defemph{\ac\ \str s}. Almost-complex \str s are also induced from \defemph{symplectic \str s}, and are themselves a particular case of a \defemph{\spinc-\str} (one can think of a \spinc-\str\ as an \ac\ \str\ defined on the $2$-skeleton).
S\pinc-\str s are used to build the \defemph{\SW\ invariants}: one can think of the latter as invariants of various \spinc-\str s. 
In the case of a symplectic manifold, the \SW\ invariant can be interpreted as a count of holomorphic curves.
The existence of symplectic \str s is essentially equivalent to the existence of \defemph{Lefschetz pencils} on the manifold. The latter can be thought of as being singular \defemph{taut foliations}.
In general, \ac\ \str s can be deformed to admit singular foliations, see \Secref{sec-sing}.
\SW\ theory yields bounds on the \defemph{minimal genus} of an embedded surface in a $4$-manifold. 
In certain special cases, taut foliations also offer genus bounds, and it is conjectured that that might be true in general.

\section{Generalities}

Let $\F$ be a $2$-dimensional foliation on a $4$-manifold $M$. It has a tangent bundle $\T{\F}$ (tangent to the leaves), and a normal bundle $\N{\F}=\T{M}\big/\T{\F}$. 
An auxiliary \riem\ metric $g$ determines an embedding of the normal bundle $\N{\F}$ as a subbundle of $\T{M}$. We then have the $g$-orthogonal splitting $\T{M}=\T{\F}\oplus\N{\F}$.

In the presence of a \riem\ metric, a foliation $\F$ induces an \emph{\ac\ \str} $J=J_\F$ on $M$: essentially $J$ rotates the plane-fibers of $\T{\F}$ and $\N{\F}$ by $\pi/2$ in the right directions. 
Or: pick a orienting orthonormal frame $\{\tau_1,\tau_2,\nu_1,\nu_2\}$ in $\T{M}$ \st\ $\tau_1, \tau_2\in\T{\F}$ and $\nu_1,\nu_2\in\N{\F}$, and define the \ac\ \str\ $J$ by $J\tau_1=\tau_2$ and $J\nu_1=\nu_2$. The leaves of $\F$ are then $J$-holomorphic, and
the homotopy class of $J_\F$ is independent of $g$.

Due to this \ac\ \str\ $J_\F$, we have well-defined Chern classes:
$c_1(\T{M},J)=c_1(\T{\F})+c_1(\N{\F})$ and 
$c_2(\T{M},J)=c_1(\T{\F})\cdot c_1(\N{\F})$.
That can be written:
\begin{gather*}
c_1(\T{M},J)=e(\T{\F})+e(\N{\F})\\
\chi(M)=e(\T{\F})\cdot e(\N{\F})
\end{gather*}
The class $c_1(\T{M})$ also has the \emph{a priori} properties of any Chern class:
\begin{gather*}
c_1(\T{M})\equiv w_2(M)\pmod{2}\\
p_1(M)=c_1(\T{M})^2-2\chi(M)
\end{gather*}
(Note that 
$p_1(M)=3\sigma(M)$, where $\sigma(M)$ is the signature of $M$.)

\section{Thurston's $h$-principle}

W.~Thurston proved in \cite{thurston} that any $2$-plane field on a $4$-mani\-fold can be homotoped to be tangent to a foliation. Namely:

\begin{theorem}[W.\ Thurston 1974]
\label{thm-thurston}
Let $\tau$ be a $2$-plane field on a manifold $M$ of dimension at least $4$. Let $K$ be a compact subset of $M$ \st\ $\tau$ is completely integrable in a neighborhood of $K$ ($K$ can be empty). Then $\tau$ is homotopic $\operatorname{rel} K$ to a completely integrable plane field.
\end{theorem} 

This theorem is also true on $3$-manifolds (see \cite{thurston-codim1}), but \emph{not} in a relative version.

The theorem is proved by first lifting the plane field to a Haefliger structure, and then deforming the latter to become a foliation using the main theorem of \cite{thurston}. The latter result has an alternate proof in \cite{eliashberg}.

Of course, as a precondition of the existence of any foliation on $M$, the tangent bundle of $M$ must split in a sum of two plane bundles (candidates for $\T{M}=\T{\F}\oplus\N{\F}$). In particular, $M$ must admit \ac\ \str s. Nonetheless, if we admit foliations with singularities, these restrictions can be circumvented.

\begin{remark}
By construction, the foliations yielded by Thurston's theorem are not taut: they contain Reeb components.
Flexibility comes with a price.
\end{remark}

\section{Taut foliations}

(Most statements of this section admit obvious generalizations to arbitrary dimensions and codimensions.)

\subsection{Definition}
A foliation $\F$ on a \riem\ manifold $(M,g)$ is called \defemph{minimal} if all its leaves are minimal surfaces in $(M,g)$ (\ie they locally minimize area; for any compact piece $K$ of a leaf, any small perturbation of $K$ rel $\del K$ will have bigger area; that is equivalent to each leaf having zero mean curvature).

A foliation $\F$ on $M$ is called \defemph{taut} if there is a \riem\ metric $g$ \st\ $\F$ is minimal in $(M,g)$  

\begin{remark}
In the special case of a codimension-$1$ foliation, tautness is equivalent to the existence of a $1$-manifold transverse to $\F$ and crossing all the leaves. A similar condition is much too strong for higher codimensions.
\end{remark}

\subsection{Mean curvature and characteristic form}

We define the \defemph{mean curvature} vector field $H\sub{\F,g}$ as (half) the normal component of $\nabla^g_{\tau_1}\tau_1+\nabla^g_{\tau_2}\tau_2$, for any orthonormal orienting frame $\{\tau_1, \tau_2\}$ in $\T{\F}$, with $\nabla^g$ the \LC\ connection of $g$. Perturbing a leaf in the direction opposite to $H$ most decreases its area. If $H=0$ along a leaf, then that leaf is a  minimal surface in $(M,g)$.

\begin{proposition}
\label{prop-char}
Let $\F$ be a foliation on $(M,g)$. Let $\mu$ be any $2$-form on $M$ \st\ $\mu\rest{\text{Leaf}}=vol_g\rest{\text{Leaf}}$ and $\mu(\tau,\nu)=0$ for any $\tau\in\T{\F}$ and $\nu\in\N{\F}$. Denote by $H$ the mean curvature vector field of the leaves.
Then
\[ d\mu(\tau_1,\tau_2,\nu)=\inner{H,\nu} \]
for any orthonormal orienting basis $\tau_1,\tau_2$ in $\T{\F}$ and any $\nu\in\N{\F}$.

In particular, $\F$ is minimal \Iff\ $d\mu(\tau_1,\tau_2,\nu)=0$ for all $\tau_1, \tau_2\in\T{\F}$ and $\nu\in\N{\F}$.
\noproof
\end{proposition}

Forms $\mu$ as above always exist, and are not unique. A particular such $\mu$ is the \defemph{characteristic form} 
$\chi=\chi\sub{\F,g}$, defined to be the unique $2$-form \st\ 
$\chi\rest{\text{Leaf}}=vol_g\rest{\text{Leaf}}$ and $\chi(x,\nu)=0$ for any $x\in\T{M}$ and $\nu\in\N{\F}$. (In other words, $\chi$ is just the volume form along the leaves.)

The condition $d\mu(\tau_1,\tau_2,\nu)=0$ is strictly weaker than asking that $\mu$ be closed. For example, for the characteristic form:
\begin{lemma}
\label{lemma-char}
If $d\chi\sub{\F,g}=0$, then $\F$ is minimal and $\F^\perp$ is integrable.
\noproof
\end{lemma}
The integrability of $\F^\perp$ is measured by the other significant component of $d\chi$, namely $d\chi(\tau,\nu_1,\nu_2)$ for $\tau\in\T{\F}$ and $\nu_1, \nu_2\in\N{\F}$. More precisely, we have
\[ d\chi(\tau,\nu_1,\nu_2)=
	\chi\bigl(\tau,[\nu_1,\nu_2]\bigr) \]

\subsection{Rummler's criterion}

In the direction opposite to the above \Propref{prop-char}, it was proved in \cite{rummler}:

\begin{theorem}[H.~Rummler 1979]
\label{thm-rummler}
Let $\F$ be a foliation on $M$. Assume that there is a $2$-form $\mu$ \st\ $\mu\rest{\text{Leaf}}>0$ and $d\mu\rest{\text{Leaf}}=0$. Then $\F$ is taut.
\end{theorem}

\noindent
When writing ``$\mu\rest{\text{Leaf}}>0$ and $d\mu\rest{\text{Leaf}}=0$'' we mean ``$\mu(\tau_1,\tau_2)>0$ and $d\mu(\tau_1,\tau_2,z)=0$ for all orienting pairs $\tau_1, \tau_2\in\T{\F}$ and $z\in\T{M}$''.

\begin{proof}
If $\mu$ is positive on the leaves, then there is a metric $g^\tau$, defined only on $\T{\F}$ \st\ $\mu\rest{\text{Leaf}}=vol_{g^\tau}\rest{\text{Leaf}}$. We need to extend $g^\tau$ to a genuine metric $g$ on $M$ that will determine a normal bundle $\N{\F}$ for $\F$ \st\ $\mu(\tau,\nu)=0$ for any $\tau\in\T{\F}$ and $\nu\in\N{\F}$. We will first find $\N{\F}$, then extend $g^\tau$. Concretely, if $\tau_1, \tau_2$ is any basis in $\T{\F}$, we define $\N{\F}$ as $\ker\mu(\tau_1,\cdot)\cap\ker\mu(\tau_2,\cdot)$. Then, after choosing some random metric $g^\nu$ on $\N{\F}$, we can extend $g^\tau$ to the whole $\T{M}=\T{\F}\oplus\N{\F}$ simply as $g=g^\tau\oplus g^\nu$. 
In view of \Propref{prop-char}, this $g$ will make $\F$ minimal.
\end{proof}

\begin{remark}
A condition stronger than tautness would be the existence of a $2$-form $\mu$ \st\ $\mu\rest{\text{Leaf}}>0$ and $d\mu=0$. The foliation $\F$ is then \defemph{calibrated} by $\mu$ (see \cite{harveylawson}), and, for a suitable metric $g$, the leaves are not only locally area minimizing, but are locally area minimizing in their homology classes (\ie if $K$ is any compact piece of a leaf and $K'$ is any immersed surface with $\del K'=\del K$ and representing the same class as $K$ in $H_2(M,\del K)$, then $K'$ has bigger area than $K$).
\end{remark}

\begin{example}
\label{ex-reeb}
\emph{A Reeb component is not taut.} 
Denote by $C$ a solid torus foliated so that the boundary torus $\del C$ is a leaf, and the interior is foliated by planes (see \Figref{fig-reeb}). Consider any $2$-dimensional foliation $\F$ (of any codimension) containing such a $C$. Assume that $\F$ is taut. Then there is a $2$-form $\mu$ satisfying the conditions from Rummler's \Thmref{thm-rummler}. But then
\[ 0=\int_C d\mu=\int_{\del C}\mu>0 \]
which is impossible.
\end{example}

\medskip

The strong link between $2$-forms and minimality of foliations is also suggested (via \ac\ \str s) by the following formula:

\begin{lemma}
\label{lemma-domega}
Let $g(x,y)=\inner{x,y}$ be a \riem\ metric on $M^{4}$ and $\nabla$ 
its \LC\ connection. Let $J$ be any $g$-orthogonal \ac\ \str, and
let $\omega(x,y)=\inner{Jx,\,y}$ be its fundamental $2$-form.
Let $x$, $z$ be any vector fields on $M$. Then:
\[ d\omega(x,Jx,z)
 =\Inner{[x,Jx],\,Jz}
 -\Inner{\Nabla{x}x+\Nabla{Jx}Jx,\ z} \]
\end{lemma}

The term $[x,Jx]$ measures the integrability of the $J$-holomorphic plane field $\aR\inner{x,Jx}$, while the normal component of the term $\Nabla{x}x+\Nabla{Jx}Jx$ is the mean curvature of the plane field $\aR\inner{x,Jx}$, and thus measures its $g$-minimality. (\Lemmaref{lemma-domega} will be proved at the end of the paper, in \Secref{sec-appendix}.)

\subsection{Sullivan's criterion}

In what follows, we wish to at least suggest the nature of Sullivan's criterion for tautness. Rigor must suffer. In a way, Sullivan's criterion  is a far-reaching generalization of the Reeb \Exref{ex-reeb} above.

A $2$-chain $\sigma$ in $M$ is called a \defemph{tangent homology} in $(M,\F)$ if there is a $3$-chain $\xi$ \st\ $\del\xi=\sigma$ and $\T{\F}\rest{\xi}\subset\T{\xi}$. (In other words, $\sigma$ is a boundary that can be ``filled-in'' using a chain $\xi$ that is tangent to $\F$.)

\begin{lemma}
If a sequence of tangent homologies converges to an (oriented) union of compact leaves of $\F$, then $\F$ is not taut.
\end{lemma}
\begin{proof}
Let $L_1,\ldots,L_k$ be some compact leaves of $\F$, and let $\sigma_n=\del\xi_n$ be a sequence of tangent homologies that converges to $L_1\cup\ldots\cup L_k$. 
Assume that $\F$ is taut, and thus there is a $2$-form $\mu$ with $\mu\rest{\text{Leaf}}>0$ and $d\mu\rest{\text{Leaf}}=0$. Then:
\[ 0=\int_{\xi_n}d\mu
	=\int_{\del\xi_n}\mu=\int_{\sigma_n}\mu
  \ \stackrel{n}{\longto}\ 
  \int_{L_1\cup\ldots\cup L_k}\mu \ >0 \]
which is impossible
\end{proof}

\noindent
(The appropriate notion of convergence here is, of course, a convergence ``in measure''.)

Of course, in general, $\F$ might not have compact leaves at all. Nonetheless, a set of non-compact leaves might come back on itself in a regular enough fashion that it could be called a ``foliation cycle''. Somewhat more rigorously, a \defemph{foliation cycle} can be defined as a transverse measure that is holonomy-invariant. Think of this measure as a way of weighting the leaves of $\F$.
Then:

\begin{theorem}[D.\ Sullivan 1979]
A foliation $\F$ is taut \Iff\ no non-trivial foliation cycle can be approximated by tangent homologies.
\noproof
\end{theorem}

This result was proved in \cite{sullivan-taut}, based on the theory developed in \cite{sullivan-cycles}. A complete recent discussion can be found in \cite[Ch.~10]{candelconlon}.
The rigorous development of these ideas involves currents.

A (rarely applicable) consequence of this criterion is:

\begin{corollary}
Let $\F$ be a foliation on $M$. If for every leaf $L$ of $\F$ there is a surface transverse to $\F$ that crosses $L$, then  $\F$ is taut. 
In fact, $\F$ can be calibrated.
\end{corollary}

Surfaces transverse to foliations are a rare occurrence, though
(unlike the codimension-$1$ case, where tautness is equivalent to the existence of transverse circles).

\section{More differential geometry}

Here are a couple of general results of a differential geometric nature, from \cite{fabiano} and \cite{garcianaveira}:

\begin{proposition}%
	[B.\ Fabiano 1984]
Let $\F$ be a minimal co\-di\-men\-sion-$2$ foliation on the \riem\ manifold $M^m$. If $M$ has non-negative Ricci curvature and $\F$ admits a complementary foliation $\F^\perp$, then: $\F$ is totally geodesic, or $e(\F^\perp)\neq 0$
\end{proposition}

\begin{proposition}%
	[F.J.\ Garc\'\i a \& A.M.\ Naveira 1995]
Let $\F$ be a minimal foliation on the \riem\ manifold $M^m$. If $M$ has non-positive Ricci curvature and $\F$ admits a complementary foliation $\F^\perp$, then: $\F$ is totally geodesic, or the leaves of $\F^\perp$ have non-zero scalar curvature.
\end{proposition}

As one can notice, these general results have pretty restrictive conditions (mainly: the existence of a complementary foliation).  In general, it is worth noticing that foliations of codimension $\geq 2$ are much harder to study that those of codimension $1$. This, in part, explains the success of foliation theory for investigating $3$-manifolds, and the lack of much of a theory for $4$-manifolds.

\section{Minimal genus of embedded surfaces}
\label{sec-min.genus}

Given any class $a\in H_2(M;\Zi)$, there are always embedded surfaces in $M$ that represent it. An open problem is to determine how simple such surfaces can be, or, in other words, what is the minimal genus that a surface representing $a$ can have. (Remember that $\chi(S)=2-2g(S)$, so minimum genus is maximum Euler characteristic.)

\subsection{Adjunctions}

Notice that, if $S$ is a $J$-holomorphic surface for some \ac\ \str\ $J$, then 
\[ \chi(S)+S\cdot S=c_1(J)\cdot S \]
(simply because $c_1(J)\cdot S=c_1(\T{M}\rest{S})=c_1(\T{S})+c_1(\N{S})=\chi(S)+S\cdot S$). 
This equality is known as the \emph{Adjunction formula} for $S$.

In general, the main and most powerful tool for obtaining genus bounds comes from \SW\ theory \cite{kronheimer.mrowka.proj,ozsvath.szabo.sympl.thom}:

\begin{proposition}[Adjunction Inequality]
Let $S$ be any embedded surface in $M$. Assume that either $M$ is of \SW\ simple type and $S$ has no sphere components, or that $S\cdot S\geq 0$.
Then, for any \SW\ basic class $\boldsymbol{\epsilon}$, we have:
\[ \chi(S)+S\cdot S \leq 
	\boldsymbol{\epsilon}\cdot S \]
In particular, if $J$ is an \ac\ \str\ admitting a symplectic structure, then 
\[ \chi(S)+S\cdot S \leq c_1(J)\cdot S \] 
\end{proposition}

Notice that one can reverse the orientation of $S$ and change only the sign of $\boldsymbol{\epsilon}\cdot S$. Thus one often writes the adjunction inequality more sharply as
\[ \chi(S)+S\cdot S \leq 
	-\module{\boldsymbol{\epsilon}\cdot S} \]

\subsection{Enter foliation}

Nonetheless, the bounds offered by \SW\ basic classes are not always sharp.  
For example, in the case of manifolds $M=N^3\times\Sph{1}$, P.~Kronheimer has proved in \cite{kronheimer} that foliations give \emph{better} bounds:

\begin{theorem}[P.\ Kronheimer 1999]
\label{thm-kron}
Consider $M=N^3\times\Sph{1}$, with $N$ a closed irreducible $3$-manifold.
Let $\Bar{\F}$ be a taut foliation in $N$, with Euler class $\bar{\boldsymbol{\epsilon}}=e(\T{\Bar{\F}})$. Let $\boldsymbol{\epsilon}$ be the image of $\bar{\boldsymbol{\epsilon}}$ in $H^2(N\times\Sph{1})$. Then, for any embedded surface $S$ in $M$, without sphere components, we have
\[ \chi(S)+S\cdot S \leq 
	\boldsymbol{\epsilon}\cdot S \]
\end{theorem}
In general, $\boldsymbol{\epsilon}$ is \emph{not} a \SW\ basic class. (Nonetheless, the proof of \Thmref{thm-kron} does use \SW\ theory: the taut foliation $\Bar{\F}$ is perturbed to a tight contact structure, which is then symplectically filled in a suitable way, and a version of the \SW\ invariants is used.)

A taut foliation $\Bar{\F}$ on $N^3$ induces an obvious product foliation $\F=\Bar{\F}\times\Sph{1}$ on $M=N\times\Sph{1}$ (with leaves $\text{Leaf}\times\{pt\}$). Then $\F$ is also  taut (pick a product metric), $\boldsymbol{\epsilon}=e(\T{\F})$ is the pull-back of $\bar{\boldsymbol{\epsilon}}=e(\T{\Bar{\F}})$, and
the \ac\ \str\ that $\F$ determines has $c_1(J_\F)=\boldsymbol{\epsilon}=e(\T{\F})$. One can then try to generalize \Thmref{thm-kron} as

\begin{conjecture}[P.\ Kronheimer 1999]
\label{conj}
Let $\F$ be a taut foliation on $M^4$, and $J_\F$ be an \ac\ \str\ induced by $\F$. Then, for any embedded surface $S$ without sphere components, we have
\[ \chi(S)+S\cdot S\leq c_1(J_\F)\cdot S \]
\end{conjecture}

\noindent
A few extra requirements are needed, \eg to exclude manifolds like $S\times\Sph{2}$. Kronheimer also proposes that the foliations be allowed singularities.

\medskip

\begin{remark}
The situation in \Thmref{thm-kron} has another peculiarity: $\F$ admits  transverse foliations. Indeed, since $\Bar{\F}$ has codimension $1$, any nowhere-zero vector-field in $N^3$ normal to $\Bar{\F}$ integrates to a $1$-dimensional foliation of $N$ that is transverse to $\Bar{\F}$. By multiplying its leaves by $\Sph{1}$, this $1$-dimensional foliation induces a $2$-dimensional foliation in $M$ that is transverse to $\F$.
\end{remark}

One could thus think of strengthening the hypothesis of \Conjref{conj} by requiring $\F$ not only to be taut, but also to admit a transverse foliation.
One might push things even further and ask that the second foliation be taut as well. But then one almost runs into:

\begin{proposition}[Two taut makes one symplectic]
\label{prop-2taut.sympl}
Let $\F$ and $\G$ be transverse foliations on $M^4$. If there is a metric $g$ that makes both $\F$ and $\G$ be minimal and orthogonal, then $M$ must admit a symplectic \str. 
Therefore, for any embedded surface $S$ we have
\[ \chi(S)+S\cdot S\leq c_1(J_\F)\cdot S \]
\end{proposition}

\noindent
(This is an immediate consequence of \Lemmaref{lemma-domega}.)

\begin{remark}
\emph{No taut, no symplectic.}
At the other extreme, if $M$ admits a \emph{non-taut} foliation $\F$, then \emph{no} \ac\ \str\ compatible with $\F$ admits symplectic \str s (or, more directly, $M$ admits no symplectic structures making the leaves of $\F$ symplectic submanifolds).
\end{remark}

\subsection{Contrast with dimension $3$}

Taut foliations on $3$-manifolds are well-understood, and are stron\-gly related to minimal genus surfaces \cite{thurston-norm,gabai}:

\begin{proposition}[W.\ Thurston 1986]
Let $N^3$ be a closed irreducible $3$-manifold, and $\F$ a taut foliation on $N$. 
For any embedded surface $S$ without sphere-components, we have
\[ \chi(S)\leq e(\T{\F})\cdot S \]
\end{proposition}

\noindent
(Usually written more sharply as $\chi(S)\leq -\module{e(\T{F})\cdot S}$.)

\begin{proposition}[D.\ Gabai 1983]
Let $N^3$ be a closed irreducible $3$-manifold. 
An embedded surface $S$ has minimal genus \Iff\ it is the leaf of a taut foliation of $N$.
\noproof
\end{proposition}

No similar statement is known for $4$-manifolds.
For a survey of foliations on $3$-manifolds, see for example \cite{gabai.survey}.

\section{Foliations and Gromov compactness}

\comment{
\begin{figure}[bth]
\begin{picture}(340,370)(0,40)

\put(0,300){\line(1,0){140}}
\put(0,305){\line(1,0){140}}
\put(0,310){\line(1,0){140}}
\put(0,315){\line(1,0){140}}
\put(0,320){\line(1,0){140}}
\put(140,300){\line(1,2){30}}
\put(140,305){\line(1,2){30}}
\put(140,310){\line(1,2){30}}
\put(140,315){\line(1,2){30}}
\put(140,320){\line(1,2){30}}
\put(0,320){\line(1,2){30}}

\put(30,380){\line(1,0){140}}

\qbezier[30](85,350)(85,380)(85,410) 
\qbezier[15](85,300)(85,285)(85,270)

\put(85,255){\makebox(0,0){$t=\pm 1$}} 

\put(170,300){\line(1,0){140}} 
\put(170,305){\line(1,0){140}}
\put(170,310){\line(1,0){140}}
\put(170,315){\line(1,0){140}}
\put(170,320){\line(1,0){140}}
\put(310,300){\line(1,2){30}}
\put(310,305){\line(1,2){30}}
\put(310,310){\line(1,2){30}}
\put(310,315){\line(1,2){30}}
\put(310,320){\line(1,2){30}}
\put(170,320){\line(1,2){30}}

\put(200,380){\line(1,0){30}} 
\put(340,380){\line(-1,0){60}}
\qbezier(200,355)(230,355)(230,385)
\qbezier(310,355)(280,355)(280,385)

\qbezier(230,385)(233,400)(255,400) 
\qbezier(280,385)(277,400)(255,400)

\qbezier[7](255,396)(255,403)(255,410) 
\qbezier[15](255,300)(255,285)(255,270)

\put(255,255){\makebox(0,0){$t=\pm 0.6$}}       

\put(0,100){\line(1,0){140}} 
\put(0,105){\line(1,0){140}}
\put(0,110){\line(1,0){140}}
\put(0,115){\line(1,0){140}}
\put(0,120){\line(1,0){140}}
\put(140,100){\line(1,2){30}}
\put(140,105){\line(1,2){30}}
\put(140,110){\line(1,2){30}}
\put(140,115){\line(1,2){30}}
\put(140,120){\line(1,2){30}}
\put(0,120){\line(1,2){30}}

\put(30,180){\line(1,0){30}} 
\put(170,180){\line(-1,0){60}}
\qbezier(30,155)(60,155)(60,185)
\qbezier(140,155)(110,155)(110,185)

\put(60,185){\line(0,1){35}} 
\put(110,185){\line(0,1){35}}

\qbezier[20](40,220)(50,220)(60,220) 
\qbezier[20](110,220)(120,220)(130,220)
\qbezier[27](60,220)(65,215)(85,215)
\qbezier[27](85,215)(105,215)(110,220)

\put(62,65){\line(0,1){15}} 
\put(108,65){\line(0,1){15}}
\qbezier(62,80)(65,95)(85,95)
\qbezier(108,80)(105,95)(85,95)

\qbezier[22](40,65)(50,65)(62,65) 
\qbezier[22](108,65)(120,65)(130,65)
\qbezier[24](62,65)(67,60)(85,60)
\qbezier[24](85,60)(103,60)(108,65)

\qbezier[5](85,90)(85,95)(85,100) 

\put(85,45){\makebox(0,0){$t=\pm 0.3$}}    

\put(170,100){\line(1,0){140}} 
\put(170,105){\line(1,0){140}}
\put(170,110){\line(1,0){140}}
\put(170,115){\line(1,0){140}}
\put(170,120){\line(1,0){140}}
\put(310,100){\line(1,2){30}}
\put(310,105){\line(1,2){30}}
\put(310,110){\line(1,2){30}}
\put(310,115){\line(1,2){30}}
\put(310,120){\line(1,2){30}}
\put(170,120){\line(1,2){30}}

\put(200,180){\line(1,0){30}} 
\put(340,180){\line(-1,0){60}}
\qbezier(200,155)(230,155)(230,185)
\qbezier(310,155)(280,155)(280,185)

\put(230,185){\line(0,1){13}} 
\put(280,185){\line(0,1){13}}
\qbezier(230,198)(230,188)(255,188)
\qbezier(280,198)(280,188)(255,188)


\qbezier(232,200)(231,190)(255,190) 
\qbezier(278,200)(279,190)(255,190)

\qbezier(234,202)(232,192)(255,192) 
\qbezier(276,202)(280,192)(255,192)

\put(235,200){\line(0,1){10}} 
\put(275,200){\line(0,1){10}}
\qbezier(235,210)(235,202)(255,202)
\qbezier(275,210)(275,202)(255,202)
\qbezier(235,210)(236,215)(242,216)
\qbezier(275,210)(274,215)(268,216)

\qbezier[5](235,208)(235,200)(255,200) 
\qbezier[5](275,208)(275,200)(255,200)
\qbezier[5](235,204)(235,196)(255,196)
\qbezier[5](275,204)(275,196)(255,196)
\qbezier[5](235,201)(235,193)(255,193)
\qbezier[5](275,201)(275,193)(255,193)

\qbezier[10](235,210)(235,215)(235,220) 
\qbezier[10](275,210)(275,215)(275,220)

\qbezier(240,206)(240,220)(255,220) 
\qbezier(270,206)(270,220)(255,220)

\put(235,100){\line(0,-1){15}} 
\put(275,100){\line(0,-1){15}}
\qbezier(235,85)(235,77)(255,77)
\qbezier(275,85)(275,77)(255,77)

\qbezier[5](235,87)(235,79)(255,79) 
\qbezier[5](275,87)(275,79)(255,79)
\qbezier[5](235,91)(235,83)(255,83)
\qbezier[5](275,91)(275,83)(255,83)
\qbezier[5](235,95)(235,87)(255,87)
\qbezier[5](275,95)(275,87)(255,87)
\qbezier[5](235,99)(235,91)(255,91)
\qbezier[5](275,99)(275,91)(255,91)
\qbezier[3](240,100)(240,95)(255,95)
\qbezier[3](270,100)(270,95)(255,95)
\qbezier[2](250,100)(240,99)(255,99)
\qbezier[2](260,100)(270,99)(255,99)

\qbezier[25](235,85)(235,72.5)(235,60) 
\qbezier[25](275,85)(275,72.5)(275,60)

\put(237,69){\line(0,1){10}} 
\put(273,69){\line(0,1){10}}
\qbezier(237,69)(237,64)(255,64)
\qbezier(273,69)(273,64)(255,64)

\qbezier(239,65)(239,62)(255,62) 
\qbezier(271,65)(271,62)(255,62)

\qbezier(241,62)(241,60)(255,60) 
\qbezier(269,62)(269,60)(255,60)

\put(255,45){\makebox(0,0){$t=0$}} 
\put(242,90){\vector(1,0){0}}
\qbezier(239,90)(220,90)(220,80)
\qbezier(220,80)(220,70)(200,70)
\put(152,68){{\small torus leaf}}

\end{picture}
\caption{Creating a torus leaf (3D movie)}
\label{fig-create.torus1}
\end{figure}
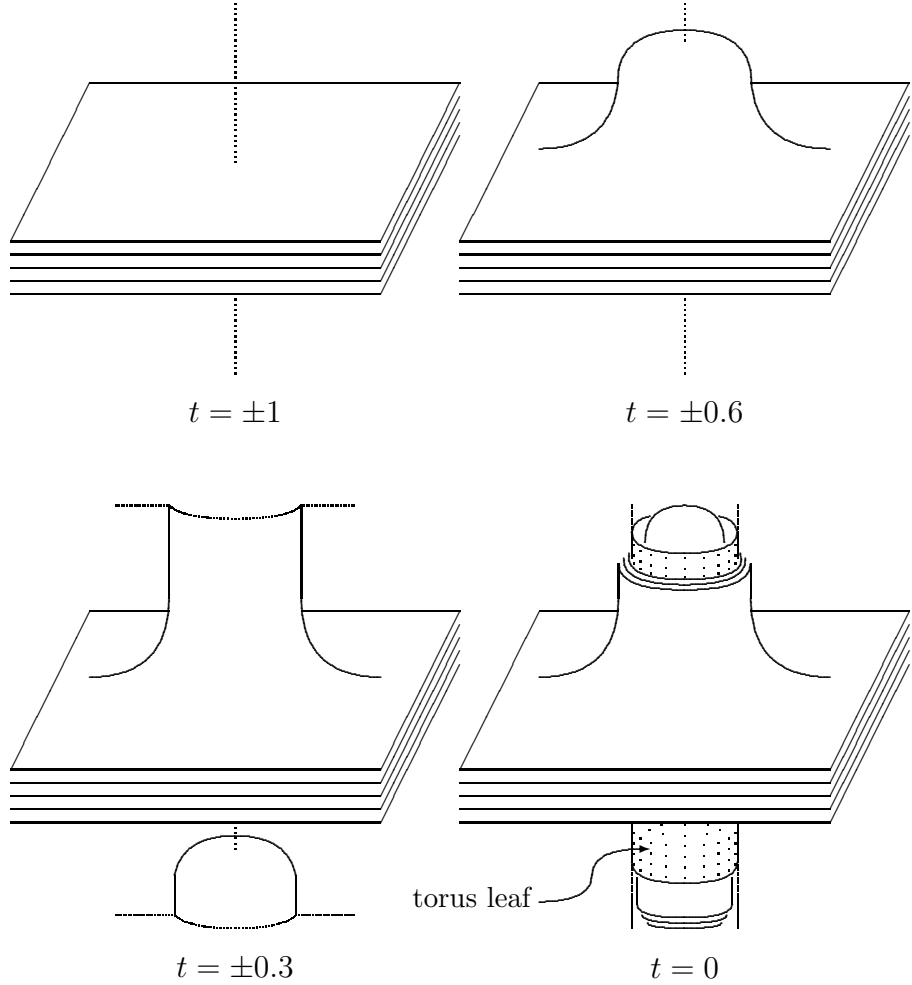
}
\comment{
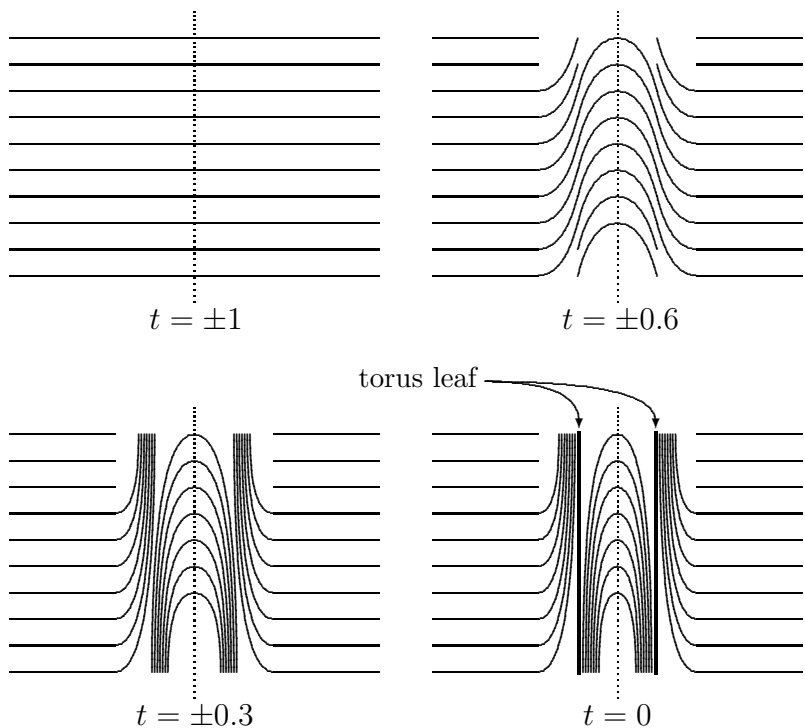
\begin{figure}[bth]
\begin{picture}(300,275)(0,25)

\put(0,200){\line(1,0){140}}
\put(0,210){\line(1,0){140}}
\put(0,220){\line(1,0){140}}
\put(0,230){\line(1,0){140}}
\put(0,240){\line(1,0){140}}
\put(0,250){\line(1,0){140}}
\put(0,260){\line(1,0){140}}
\put(0,270){\line(1,0){140}}
\put(0,280){\line(1,0){140}}
\put(0,290){\line(1,0){140}}

\qbezier[55](70,190)(70,255)(70,300)

\put(52,180){$t=\pm 1$}

\put(160,200){\line(1,0){40}}
\put(160,210){\line(1,0){40}}
\put(160,220){\line(1,0){40}}
\put(160,230){\line(1,0){40}}
\put(160,240){\line(1,0){40}}
\put(160,250){\line(1,0){40}}
\put(160,260){\line(1,0){40}}
\put(160,270){\line(1,0){40}}
\put(160,280){\line(1,0){40}}
\put(160,290){\line(1,0){40}}

\put(300,200){\line(-1,0){40}}
\put(300,210){\line(-1,0){40}}
\put(300,220){\line(-1,0){40}}
\put(300,230){\line(-1,0){40}}
\put(300,240){\line(-1,0){40}}
\put(300,250){\line(-1,0){40}}
\put(300,260){\line(-1,0){40}}
\put(300,270){\line(-1,0){40}}
\put(300,280){\line(-1,0){40}}
\put(300,290){\line(-1,0){40}}

\qbezier(200,270)(210,270)(215,290)
\qbezier(245,290)(250,270)(260,270)

\qbezier(200,260)(210,260)(215,280)
\qbezier(245,280)(250,260)(260,260)

\qbezier(200,250)(210,250)(215,270)
\qbezier(215,270)(220,290)(230,290)
\qbezier(230,290)(240,290)(245,270)
\qbezier(245,270)(250,250)(260,250)

\qbezier(200,240)(210,240)(215,260)
\qbezier(215,260)(220,280)(230,280)
\qbezier(230,280)(240,280)(245,260)
\qbezier(245,260)(250,240)(260,240)

\qbezier(200,230)(210,230)(215,250)
\qbezier(215,250)(220,270)(230,270)
\qbezier(230,270)(240,270)(245,250)
\qbezier(245,250)(250,230)(260,230)

\qbezier(200,220)(210,220)(215,240)
\qbezier(215,240)(220,260)(230,260)
\qbezier(230,260)(240,260)(245,240)
\qbezier(245,240)(250,220)(260,220)

\qbezier(200,210)(210,210)(215,230)
\qbezier(215,230)(220,250)(230,250)
\qbezier(230,250)(240,250)(245,230)
\qbezier(245,230)(250,210)(260,210)

\qbezier(200,200)(210,200)(215,220)
\qbezier(215,220)(220,240)(230,240)
\qbezier(230,240)(240,240)(245,220)
\qbezier(245,220)(250,200)(260,200)

\qbezier(215,210)(220,230)(230,230)
\qbezier(230,230)(240,230)(245,210)

\qbezier(215,200)(220,220)(230,220)
\qbezier(230,220)(240,220)(245,200)

\qbezier[55](230,190)(230,255)(230,300)

\put(208,180){$t=\pm 0.6$}

\put(0,50){\line(1,0){40}}
\put(0,60){\line(1,0){40}}
\put(0,70){\line(1,0){40}}
\put(0,80){\line(1,0){40}}
\put(0,90){\line(1,0){40}}
\put(0,100){\line(1,0){40}}
\put(0,110){\line(1,0){40}}
\put(0,120){\line(1,0){40}}
\put(0,130){\line(1,0){40}}
\put(0,140){\line(1,0){40}}

\put(140,50){\line(-1,0){40}}
\put(140,60){\line(-1,0){40}}
\put(140,70){\line(-1,0){40}}
\put(140,80){\line(-1,0){40}}
\put(140,90){\line(-1,0){40}}
\put(140,100){\line(-1,0){40}}
\put(140,110){\line(-1,0){40}}
\put(140,120){\line(-1,0){40}}
\put(140,130){\line(-1,0){40}}
\put(140,140){\line(-1,0){40}}

\qbezier(40,50)(55,50)(55,140)
\qbezier(40,60)(54,60)(54,140)
\qbezier(40,70)(53,70)(53,140)
\qbezier(40,80)(52,80)(52,140)
\qbezier(40,90)(51,90)(51,140)
\qbezier(40,100)(50,100)(50,140)
\qbezier(40,110)(49,110)(49,140)
 
\qbezier(100,50)(85,50)(85,140)
\qbezier(100,60)(86,60)(86,140)
\qbezier(100,70)(87,70)(87,140)
\qbezier(100,80)(88,80)(88,140)
\qbezier(100,90)(89,90)(89,140)
\qbezier(100,100)(90,100)(90,140)
\qbezier(100,110)(91,110)(91,140)

\qbezier(54,50)(54,140)(70,140)
\qbezier(55,50)(55,130)(70,130)
\qbezier(56,50)(56,120)(70,120)
\qbezier(57,50)(57,110)(70,110)
\qbezier(58,50)(58,100)(70,100)
\qbezier(59,50)(59,90)(70,90)
\qbezier(60,50)(60,80)(70,80)

\qbezier(86,50)(86,140)(70,140)
\qbezier(85,50)(85,130)(70,130)
\qbezier(84,50)(84,120)(70,120)
\qbezier(83,50)(83,110)(70,110)
\qbezier(82,50)(82,100)(70,100)
\qbezier(81,50)(81,90)(70,90)
\qbezier(80,50)(80,80)(70,80)

\qbezier[55](70,40)(70,105)(70,150)

\put(47,30){$t=\pm 0.3$}

\put(160,50){\line(1,0){40}}
\put(160,60){\line(1,0){40}}
\put(160,70){\line(1,0){40}}
\put(160,80){\line(1,0){40}}
\put(160,90){\line(1,0){40}}
\put(160,100){\line(1,0){40}}
\put(160,110){\line(1,0){40}}
\put(160,120){\line(1,0){40}}
\put(160,130){\line(1,0){40}}
\put(160,140){\line(1,0){40}}

\put(300,50){\line(-1,0){40}}
\put(300,60){\line(-1,0){40}}
\put(300,70){\line(-1,0){40}}
\put(300,80){\line(-1,0){40}}
\put(300,90){\line(-1,0){40}}
\put(300,100){\line(-1,0){40}}
\put(300,110){\line(-1,0){40}}
\put(300,120){\line(-1,0){40}}
\put(300,130){\line(-1,0){40}}
\put(300,140){\line(-1,0){40}}

\qbezier(200,50)(214,50)(214,140)
\qbezier(200,60)(213,60)(213,140)
\qbezier(200,70)(212,70)(212,140)
\qbezier(200,80)(211,80)(211,140)
\qbezier(200,90)(210,90)(210,140)
\qbezier(200,100)(209,100)(209,140)
\qbezier(200,110)(208,110)(208,140)
 
\qbezier(260,50)(246,50)(246,140)
\qbezier(260,60)(247,60)(247,140)
\qbezier(260,70)(248,70)(248,140)
\qbezier(260,80)(249,80)(249,140)
\qbezier(260,90)(250,90)(250,140)
\qbezier(260,100)(251,100)(251,140)
\qbezier(260,110)(252,110)(252,140)

\qbezier(217,50)(217,140)(230,140)
\qbezier(218,50)(218,130)(230,130)
\qbezier(219,50)(219,120)(230,120)
\qbezier(220,50)(220,110)(230,110)
\qbezier(221,50)(221,100)(230,100)
\qbezier(222,50)(222,90)(230,90)
\qbezier(223,50)(223,80)(230,80)

\qbezier(243,50)(243,140)(230,140)
\qbezier(242,50)(242,130)(230,130)
\qbezier(241,50)(241,120)(230,120)
\qbezier(240,50)(240,110)(230,110)
\qbezier(239,50)(239,100)(230,100)
\qbezier(238,50)(238,90)(230,90)
\qbezier(237,50)(237,80)(230,80)

\linethickness{1pt}
\put(215.5,49){\line(0,1){92}}
\put(244.5,49){\line(0,1){92}}
\thinlines
 
\qbezier[55](230,40)(230,105)(230,150)

\qbezier(215.5,145)(215.5,160)(180,160)
\put(215.5,142){\vector(0,-1){0}}
\qbezier(244.5,145)(244.5,160)(180,160)
\put(244.5,142){\vector(0,-1){0}}
\put(132,158){{\small torus leaf}}

\put(216,30){$t=0$}

\end{picture}
\caption{Creating a torus leaf (2D movie)}
\label{fig-create.torus2}
\end{figure}
}

First, an example that shows the flexibility of foliations:

\begin{example}
\label{ex-create.torus}
\emph{Creating a torus leaf.}
Let  $\F$ be \emph{any} foliation on a $4$-manifold $M$. Let $c:\Sph{1}\to M$ be any embedding. The curve $c$ can always be slightly perturbed to be transverse to $\F$. Choose another local coordinate near $c$, transverse both to $\F$ and to $c$, and think of it as \emph{time} (with $c$ appearing at time $t=0$). Start at time $t=-1$. As time goes on, begin pushing more and more the leaves of $\F$ parallel with the direction of $c$, wrapping them around more and more as time approaches $t=0$ 
(see Figures \ref{fig-create.torus1} and \ref{fig-create.torus2}).
At $t=0$, we can fit in a torus leaf, with the interior of the torus foliated by leaves diffeomorphic to $\aR^2$---a \emph{Reeb component}. As time goes on from $t=0$, we play the movie backward. Notice that the new foliation is \emph{homotopic} with the one we started with (compare with \Corref{cor-torus}).
\end{example}

Combining with \Thmref{thm-exist} below, this yields:

\begin{corollary}
\label{cor-no.gromov}
Any \ac\ \str\ is homotopic to one for which Gromov compactness fails.
\end{corollary}

\defemph{Gromov compactness} here means the compactness of the space of all $J$-holomorphic curves (curve = real surface). In other words, any sequence of holomorphic curves $f_n:(\Sigma,j_n)\to (M,J)$ has a subsequence converging to a limit $f:\Sigma^*\to M$ that is $J$-holomorphic (and may have nodal singularities, see also \Figref{fig-bubble}; the limit domain $\Sigma^*$ is obtained by collapsing circles of $\Sigma$). 
Gromov compactness always holds for symplectic \str s, due to the following \cite{gromov.pseudo}:

\begin{figure}[bth]
\setlength{\unitlength}{0.7pt}
\begin{picture}(200,200)(-100,-100)

\thicklines
\put(-100,0){\line(1,0){200}}
\put(0,-100){\line(0,1){200}}
\thinlines

\put(-100,-7){\line(1,0){50}}
\qbezier(-50,-7)(-7,-7)(-7,-50)
\put(-7,-50){\line(0,-1){50}}
\qbezier(-100,-15)(-10,-10)(-15,-100)
\qbezier(-100,-25)(-30,-30)(-25,-100)
\qbezier(-100,-35)(-50,-50)(-35,-100)

\put(100,7){\line(-1,0){50}}
\qbezier(50,7)(7,7)(7,50)
\put(7,50){\line(0,1){50}}
\qbezier(100,15)(10,10)(15,100)
\qbezier(100,25)(30,30)(25,100)
\qbezier(100,35)(50,50)(35,100)

\put(100,-50){%
	\makebox(0,0){$z_1 z_2=\varepsilon\to 0$}}

\end{picture}
\caption{Apparition of singularities in the limit}
\label{fig-bubble}
\end{figure}
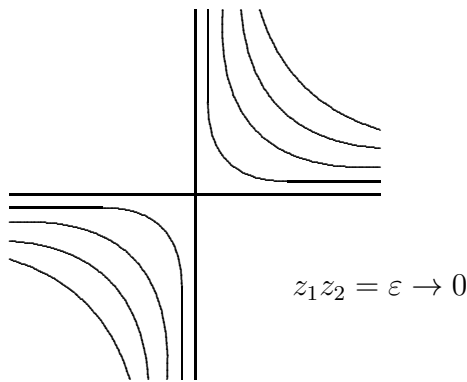

\begin{theorem}[Gromov compactness]
Let $J$ be an \ac\ \str\ on $M$, and $g$ a compatible \riem\ metric.   
Let $f_n:(\Sigma,j_n)\to (M,J)$ be a sequence
of holomorphic compact curves in $M$.
If the $g$-area of $f_n(\Sigma)$ is bounded, then
$f_n$ has a subsequence that converges to
a holomorphic curve $f:(\Sigma^*,j)\to(M,J)$, which might have nodal singularities.
\end{theorem}

\noindent
And the area bound condition is satisfied if $J$ admits a compatible \emph{symplectic} \str, and if all $f_n(\Sigma)$ represent a same homology class.
(This is actually how the theorem is most used.)
For a thorough discussion, see for example \cite{audin.lafontaine}.

\begin{Proof}{\Corref{cor-no.gromov}}
By \Thmref{thm-exist}, an \ac\ \str\ can be deformed till there is a singular foliation $\F$ with all leaves holomorphic. 
As in \Exref{ex-create.torus} above, create a torus leaf. 
Actually, by ``freezing'' the movie at $t=0$ (expanding the frame at $t=0$ to all $t\in[-\epsilon,\epsilon]$), create a lot of tori.
Now pick a second curve, orthogonal to these tori, and apply that example again.  What appears in the end is a torus that explodes to make room for a new Reeb component. Thinking in terms of an \ac\ \str\ $J$ induced by the final foliation, we have a sequence of $J$-holomorphic tori that has no decent limit.
\end{Proof}

\begin{question}[R.~Kirby 2002]
\label{q-kirby}
Let $\F$ be a foliation on $M$, and $J$ an \ac\ \str\ making the leaves $J$-holomorphic. What conditions imposed on $\F$ insure that Gromov compactness holds for $J$?
\end{question}

In the extreme, if $\F$ is a Lefschetz pencil (see \Exref{ex-donaldson}), then Gromov compactness holds (the manifold is symplectic). Compare also with \Propref{prop-2taut.sympl}.

The flexibility from our examples, at least, is done away with if we require the foliations to be \emph{taut} (since that excludes Reeb components).

\section{Singular foliations}
\label{sec-sing}

For a fuller treatment of the topic of this section, see \cite{scorpan.foliation.exist}.

\comment{
\begin{figure}[hbt]
\setlength{\unitlength}{0.9pt}
\begin{picture}(300,320)(0,0)

\qbezier(150,300)(-100,200)(150,100)
\qbezier(150,300)(-80,200)(150,100)
\qbezier(150,300)(-60,200)(150,100)
\qbezier(150,300)(-40,200)(150,100)
\qbezier(150,300)(-20,200)(150,100)
\qbezier(150,300)(0,200)(150,100)
\qbezier(150,300)(20,200)(150,100)
\qbezier(150,300)(40,200)(150,100)
\qbezier(150,300)(60,200)(150,100)
\qbezier(150,300)(80,200)(150,100)
\qbezier(150,300)(100,200)(150,100)
\qbezier(150,300)(120,200)(150,100)
\qbezier(150,300)(140,200)(150,100)
\qbezier(150,300)(160,200)(150,100)
\qbezier(150,300)(180,200)(150,100)
\qbezier(150,300)(200,200)(150,100)
\qbezier(150,300)(220,200)(150,100)
\qbezier(150,300)(240,200)(150,100)
\qbezier(150,300)(260,200)(150,100)
\qbezier(150,300)(280,200)(150,100)

\qbezier(150,300)(185,280)(215,240)
\qbezier(150,100)(185,120)(215,160)

\qbezier(215,240)(260,170)(290,170)
\qbezier(215,160)(260,230)(290,230)

\qbezier(290,170)(330,167)(330,200)
\qbezier(290,230)(330,233)(330,200)

\qbezier(150,300)(190,280)(230,240)
\qbezier(150,100)(190,120)(230,160)

\qbezier(150,300)(195,280)(245,240)
\qbezier(150,100)(195,120)(245,160)

\qbezier(150,300)(200,280)(260,240)
\qbezier(150,100)(200,120)(260,160)

\put(150,50){\makebox(0,0){$\boldsymbol{\downarrow}$}}

\linethickness{2pt}
\put(30,0){\line(1,0){240}}
\put(280,-5){$\mathbb{C}P^1$}

\end{picture}
\caption{A Lefschetz pencil}
\label{fig-lefschetz}
\end{figure}
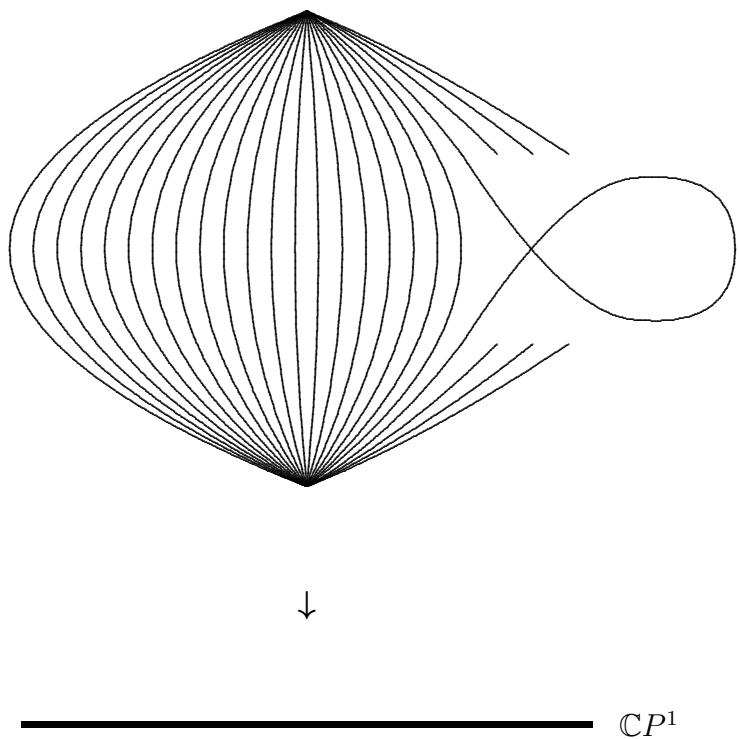
}

For a foliation $\F$ to exist on a manifold $M$, the tangent bundle must split $\T{M}=\T{\F}\oplus\N{\F}$. Since in general that does not happen, one must allow for singularities of $\F$. An important example is \cite{donaldson.lefschetz}:

\begin{theorem}[S.K.\ Donaldson]
\label{ex-donaldson}
Let $J$ be an \ac\ \str\ on $M$ that admits a compatible symplectic structure. Then $J$ can be deformed to an \ac\ \str\ $J'$ \st\ $M$ admits a Lefschetz pencil with $J'$-holomorphic fibers. 
\end{theorem}

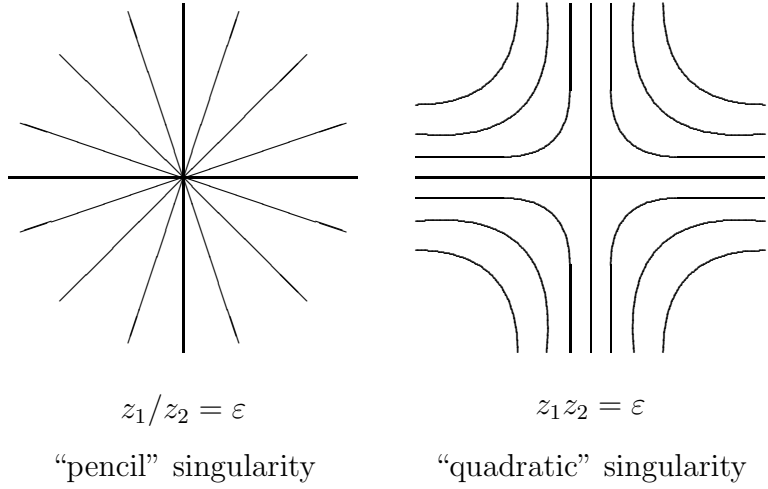
\begin{figure}[bth]
\setlength{\unitlength}{1.1pt}
\begin{picture}(260,170)(-50,-55)

\put(-50,50){\line(1,0){120}}
\put(10,-10){\line(0,1){120}}
\put(-32.5,7.5){\line(1,1){85}}
\put(-32.5,92.5){\line(1,-1){85}}

\put(10,50){\line(3,1){56}} 
\put(10,50){\line(-3,1){56}}
\put(10,50){\line(3,-1){56}}
\put(10,50){\line(-3,-1){56}}

\put(10,50){\line(1,3){19}}
\put(10,50){\line(-1,3){19}}
\put(10,50){\line(1,-3){19}}
\put(10,50){\line(-1,-3){19}}

\put(10,-30){\makebox(0,0){$z_1/z_2=\varepsilon$}}
\put(10,-50){\makebox(0,0){``pencil'' singularity}}

\put(90,50){\line(1,0){120}}
\put(150,110){\line(0,-1){120}}

\put(90,57){\line(1,0){30}}
\qbezier(120,57)(143,57)(143,80)
\put(143,80){\line(0,1){30}}

\put(90,43){\line(1,0){30}}
\qbezier(120,43)(143,43)(143,20)
\put(143,20){\line(0,-1){30}}

\put(210,57){\line(-1,0){30}}
\qbezier(180,57)(157,57)(157,80)
\put(157,80){\line(0,1){30}}

\put(210,43){\line(-1,0){30}}
\qbezier(180,43)(157,43)(157,20)
\put(157,20){\line(0,-1){30}}

\qbezier(125,110)(125,75)(90,75)
\qbezier(175,110)(175,75)(210,75)
\qbezier(125,-10)(125,25)(90,25)
\qbezier(175,-10)(175,25)(210,25)

\qbezier(135,110)(140,60)(90,65)
\qbezier(165,110)(160,60)(210,65)
\qbezier(135,-10)(140,40)(90,35)
\qbezier(165,-10)(160,40)(210,35)

\put(150,-30){\makebox(0,0){$z_1 z_2=\varepsilon$}}
\put(150,-50){\makebox(0,0){``quadratic'' singularity}}

\end{picture}
\caption{Singularities in a Lefschetz pencil}
\label{fig-sing}
\end{figure}

A \defemph{Lefschetz pencil} is a \emph{singular fibration} $M\to\CP^1$ (see Figures \ref{fig-lefschetz} and \ref{fig-sing}) with singularities modeled locally by 
\[ (z_1,z_2)\maps z_1/z_2  \qquad\text{or}\qquad (z_1,z_2)\maps z_1 z_2 \]
for suitable local complex coordinates (compatible with the orientation of $M$). Note that all fibers pass through all singularities of type $z_1/z_2$. 
The existence of a Lefschetz pencil is equivalent to the existence of a symplectic structure. See \cite[ch.~8]{gompf} for a survey.

One can think of a Lefschetz pencil as a special case of a taut \emph{singular foliation}.

\medskip

As mentioned before, a (non-singular) foliation $\F$ on $M$ induces \ac\ \str s, and we have
$c_1(J_\F)=e(\T{\F})+e(\N{\F})$ and
$\chi(M)=e(\T{\F})\cdot e(\N{\F})$
If the foliation $\F$ has singularities, then the second equality above fails, and the defect $\chi(M)-e(\T{\F})\cdot e(\N{\F})$ measures the number of singularities (or, for more general singularities, their complexity).

Call a class $c\in H^{2}(M;\Zi)$ a \defemph{complex class} of the $4$-manifold $M$ if
\[ c\equiv w_2(M) \pmod{2} \qquad\text{and}\qquad 
  p_1(M)=c^2-2\chi(M) \]
An element $c\in H^2(M;\Zi)$ is a  complex class \Iff\ there is an \ac\ \str\ $J$ on $M$ \st\ $c_1(\T{M},J)=c$. 

\begin{theorem}[Existence theorem]
\label{thm-exist}
Let $c\in H^2(M;\Zi)$ be a complex class, and let $c=\tau+\nu$ be any splitting \st\ $\chi(M)-\tau\nu\geq 0$.
Then there is a singular foliation $\F$ with $e(\T{\F})=\tau$ and $e(\N{\F})=\nu$, with $n=\chi(M)-\tau\nu$ singularities that can be chosen to be modeled on the levels of the complex functions $(z_1,z_2)\maps z_1/z_2$ or $(z_1,z_2)\maps z_1 z_2$.
\end{theorem}

\begin{remark}
Due to the singularities, the bundles $\T{\F}$ and $\N{\F}$ are only defined on $M\setminus\{\text{singularities}\}$. Their Euler classes \emph{a priori} belong to $H^2(M\setminus\{\text{singularities}\};\,\Zi)$, but can be pulled-back to $H^2(M;\Zi)$, since the isolated singularities can be chosen to affect only the 4-skeleton of $M$, and thus to not influence $H^2$.
\end{remark}

\begin{remark}
Unlike a Lefschetz pencil, in general \emph{not all} leaves of the foliation pass through the $z_1/z_2$-singularities. If they did, then the foliation would be taut, which is rather rare.
\end{remark}

\medskip

Finding a splitting $c=\tau+\nu$ with $\chi(M)-\tau\nu\geq 0$ is possible for most $4$-manifolds that admit \ac\ \str s. For example, if 
\[ \chi(M)\geq 0 \] 
(\eg for all \emph{simply-connected} $M$'s), then one can choose either one of $\tau$ or $\nu$ to be $0$, and conclude that such foliations exist. Or:

\begin{lemma}
If $b_2^+(M)>0$, then there are infinitely many splittings $c=\tau+\nu$ with $\chi(M)-\tau\nu\geq 0$ (and thus infinitely many homotopy types of foliations). If $b_2^+(M)=0$, there are at most finitely many.
\end{lemma}

\begin{proof}
If $b_2^+(M)>0$, there is a class $\alpha$ with $\alpha\cdot\alpha>0$. Choose $\tau=c-k\alpha$  and $\nu=k\alpha$ ($k\in\Zi$). 
Then $\chi(M)-\tau\nu=\chi(M)-kc\alpha+k^2\alpha^2$, and for $k$ big enough it will be positive.
\end{proof}

The main restriction to the existence of such foliations remains, of course, the existence of a complex class.
But \Thmref{thm-exist} can be generalized for the case when $c$ is merely an \emph{integral lift} of $w_2(M)$ (corresponding to weakening \ac\ \str s to \emph{\spinc-\str s}). 
In that case, singularities are also modeled using local complex coordinates, but are allowed to be compatible either with the orientation of $M$ or with the \emph{opposite orientation}.
This is similar to the generalization of Lefschetz pencils to
\emph{achiral Lefschetz pencils}, see \cite[\S8.4]{gompf}.

\medskip

The singularities of $\F$ are exactly the singularities that appear in a Lefschetz pencil. They can be chosen in either combination of types as long as their number is $n=\chi(M)-\tau\nu$. For example, there are always foliations with only $z_1/z_2$-singularities, that can thus be eliminated by blowing-up.
In fact, other choices of singularities are possible.

Namely, for any isolated singularity $p$ of a foliation that is compatible with an \ac\ \str\ one can define its \defemph{Hopf degree} $\deg p\geq 0$ (essentially a Hopf invariant of the tangent plane field above a small $3$-sphere around $p$), and
then one can choose any suitable set of singularities $\{p_i\}$ as long as 
\[ \tsum \deg p_i = \chi(M)-\tau\nu \]

\medskip

Let $\F$ be a foliation. If $S$ is a closed transversal of $\F$, then we must have 
\begin{align*}
& e(\T{\F})\cdot S=e(\T{\F}\rest{S})=e(\N{S})=S\cdot S\\
& e(\N{\F})\cdot S=e(\N{\F}\rest{S})=e(\T{S})=\chi(S)
\end{align*}
These conditions are, in fact, sufficient:

\begin{theorem}[Closed transversal]
\label{thm-transversal}
Let $S$ be a closed connected surface. 
Let $c$ be a complex class with a splitting  $c=\tau+\nu$ \st\ $\chi(M)-\tau\nu\geq 0$. 
If
\[ \chi(S)=\nu\cdot S \qquad\qquad 
  S\cdot S=\tau\cdot S \]
then there is a singular foliation $\F$ with $e(\T{\F})=\tau$, $e(\N{\F})=\nu$, and having $S$ as a closed transversal.
\end{theorem}

If, on the other hand, $S$ is a closed leaf of $\F$, then we have 
\begin{align*}
& e(\T{\F})\cdot S=e(\T{\F}\rest{S})=e(\T{S})=\chi(S)\\
& e(\N{\F})\cdot S=e(\N{\F}\rest{S})=e(\N{S})=S\cdot S
\end{align*}
Conversely:

\begin{theorem}[Closed leaf]
\label{thm-leaf}
Let $S$ be a closed connected surface with $S\cdot S\geq 0$. 
Let $c$ be a complex class with a splitting $c=\tau+\nu$ \st\  $\chi(M)-\tau\nu\geq S\cdot S$.
If
\[ \chi(S)=\tau\cdot S \qquad\qquad 
  S\cdot S=\nu\cdot S \]
then there is a singular foliation $\F$ with $e(\T{\F})=\tau$, $e(\N{\F})=\nu$, and having $S$ as a closed leaf.
(The number of singularities along $S$ is $S\cdot S$.)
\end{theorem}

An immediate consequence is:

\begin{corollary}[Trivial tori]
\label{cor-torus}
A homologically-trivial torus can always be made a leaf or a transversal of a foliation.
\end{corollary}

Such flexibility is a strong suggestion that more rigidity is needed in order to actually catch any of the topology of $M$ with the aid of foliations. Requiring foliations to be taut seems, once again, a natural suggestion.

\section{Appendix}
\label{sec-appendix}

\subsection{Minimality and $2$-forms}
\begin{Proof}{\Lemmaref{lemma-domega}}
We prove that, if $g$ is a \riem\ metric, $\nabla$ 
its \LC\ connection, $J$ be any $g$-orthogonal \ac\ \str, and
$\omega(x,y)=\inner{Jx,\,y}$ its fundamental $2$-form,
then, for any vector fields $x, z$ on $M$, we have:
\[ (d\omega)(x,Jx,z)
 =\Inner{[x,Jx],\,Jz}
 -\Inner{\Nabla{x}x+\Nabla{Jx}Jx,\ z} \]

For any $2$-form $\alpha$ we have:
\begin{align*}
(d\alpha)(x,y,z)
&=\bigl(\Nabla{x}\alpha\bigr)(y,z)
 +\bigl(\Nabla{y}\alpha\bigr)(z,x)
 +\bigl(\Nabla{z}\alpha\bigr)(x,y)\\
&=x\alpha(y,z)+y\alpha(z,x)+z\alpha(x,y)\\
&\qquad
 -\alpha\bigl(\Nabla{x}y,\,z\bigr)-\alpha\bigl(\Nabla{y}z,\,x\bigr)
 -\alpha\bigl(\Nabla{z}x,\,y\bigr)\\
&\qquad
 -\alpha\bigl(y,\,\Nabla{x}z\bigr)-\alpha\bigl(z,\,\Nabla{y}x\bigr)
 -\alpha\bigl(x,\,\Nabla{z}y\bigr)
\end{align*}
Applying this to $\omega(a,b)=\inner{Ja,\,b}$, we have:
\begin{align*}
(d\omega)(x,Jx,z)
&=-x\Inner{x,z}+(Jx)\Inner{Jz,\,x}+z\Inner{Jx,Jx}\\
&\qquad
 -\Inner{J\Nabla{x}Jx,\,z}-\Inner{J\Nabla{Jx}z,\,x}
 -\Inner{J\Nabla{z}x,\,Jx}\\
&\qquad
 +\Inner{x,\,\Nabla{x}z}-\Inner{Jz,\,\Nabla{Jx}x}
 -\Inner{Jx,\,\Nabla{z}Jx}
\end{align*}
Using that $\inner{a,Jb}=-\inner{Ja,b}$ we get: 
\begin{align*}
(d\omega)(x,Jx,z)
&=-x\Inner{x,z}+(Jx)\Inner{Jz,\,x}+z\Inner{x,x}\\
&\qquad
 +\Inner{\Nabla{x}Jx,\,Jz}+\Inner{\Nabla{Jx}z,\,Jx}
 -\Inner{\Nabla{z}x,\,x}\\
&\qquad
 +\Inner{x,\,\Nabla{x}z}-\Inner{Jz,\,\Nabla{Jx}x}
 -\Inner{Jx,\,\Nabla{z}Jx}
\end{align*}
Since 
$\Inner{\Nabla{z}x,\,x}=\rec{2}z\Inner{x,x}$
and
$\Inner{Jx,\,\Nabla{z}Jx}=\rec{2}z\Inner{Jx,Jx}=\rec{2}z\Inner{x,x}$, 
we cancel the last terms of  each line, and get:
\begin{align*}
(d\omega)(x,Jx,z)
&=-x\Inner{x,z}+(Jx)\Inner{Jz,\,x}\\
&\qquad
 +\Inner{\Nabla{x}Jx,\,Jz}+\Inner{\Nabla{Jx}z,\,Jx}\\
&\qquad
 +\Inner{x,\,\Nabla{x}z}-\Inner{Jz,\,\Nabla{Jx}x}
\end{align*}
But $x\Inner{x,z}=\Inner{\Nabla{x}x,\,z}+\Inner{x,\,\Nabla{x}z}$, so
$(Jx)\Inner{Jz,\,x} = -(Jx)\Inner{z,\,Jx}
 = -\Inner{\Nabla{Jx}z,\,Jx}-\Inner{z,\,\Nabla{Jx}Jx}$, 
and therefore:
\begin{align*}
(d\omega)(x,Jx,z)
&=-\Inner{\Nabla{x}x,\,z}-\Inner{x,\,\Nabla{x}z}
 -\Inner{\Nabla{Jx}z,\,Jx}-\Inner{z,\,\Nabla{Jx}Jx}\\
&\quad
 +\Inner{\Nabla{x}Jx,\,Jz}+\Inner{\Nabla{Jx}z,\,Jx}
 +\Inner{x,\,\Nabla{x}z}-\Inner{Jz,\,\Nabla{Jx}x}\\
(d\omega)(x,Jx,z)
&=-\Inner{\Nabla{x}x,\,z}-\Inner{z,\,\Nabla{Jx}Jx}
 +\Inner{\Nabla{x}Jx,\,Jz}-\Inner{Jz,\,\Nabla{Jx}x}
\end{align*}
Since $\nabla$ is torsion-free, we have 
$\Nabla{x}Jx-\Nabla{Jx}x=[x,Jx]$, so:
\[ (d\omega)(x,Jx,z)
 =\Inner{[x,Jx],\,Jz}
 -\Inner{\Nabla{x}x+\Nabla{Jx}Jx,\ z} \]
which concludes the proof.
\end{Proof}%

In particular, if $\omega$ is symplectic (\ie $d\omega=0$), then any 
integrable $J$-holomorphic plane field is $g$-minimal, and, vice-versa, any $g$-minimal 
$J$-holomorphic plane field must be integrable. 
The converse is also true: If there are enough $J$-holomorphic integrable minimal 
plane fields, then $\omega$ must be symplectic. 
Thus:

\begin{corollary}
Assume that $M$ admits two transversal $2$-dimensional foliations 
$\F$ and $\G$ such that: 
there is a metric $g$ \st\ both $\F$ and $\G$ are 
$g$-minimal;
and there is an $g$-orthogonal \ac\ \str\ $J$ that makes both $\F$ 
and $\G$ be $J$-holomorphic.
Then $M$ admits the symplectic structure 
$\omega(x,y)=g(Jx,y)$.
\end{corollary}

In particular, if the first condition is satisfied, and further 
$\F$ and $\G$ are $g$-orthogonal, then the second condition is 
automatically satisfied, and \Propref{prop-2taut.sympl} follows:

\emph{
If a \riem\ manifold $M$ admits two $g$-orthogonal and $g$-minimal 
foliations, then $M$ admits a symplectic structure.
}

\subsection{Foliated surgery}

We add here a few comments on some surgical techniques that can be used to combine foliations.

\begin{piece}
Transverse surfaces to foliations can be used to glue foliated manifolds using a \emph{normal connected sum}:
Let $(M,\F)$ and $(N,\G)$ be foliated manifolds. Assume that there is an embedded surface $S_1$ in $M$ transverse to $\F$, and an embedded surface $S_2$ in $M$ transverse to $\G$. Assume further that $S_1\iso S_2$ and $S_1\cdot S_1=-S_2\cdot S_2$. Then the foliated manifolds $(M,\F)$ and $(N,\G)$ admit a foliated normal connected sum along $S_1\iso S_2$
\end{piece}

\begin{piece}
For regular \emph{connected sums}, one can use pencil singularities of opposite orientations. Namely, when doing the connected sum, pick small $3$-spheres around two pencil singularities of opposite types. The result of the connected sum is a foliated manifold with one less positive pencil singularity and one less negative pencil singularity.
\end{piece}

\comment{
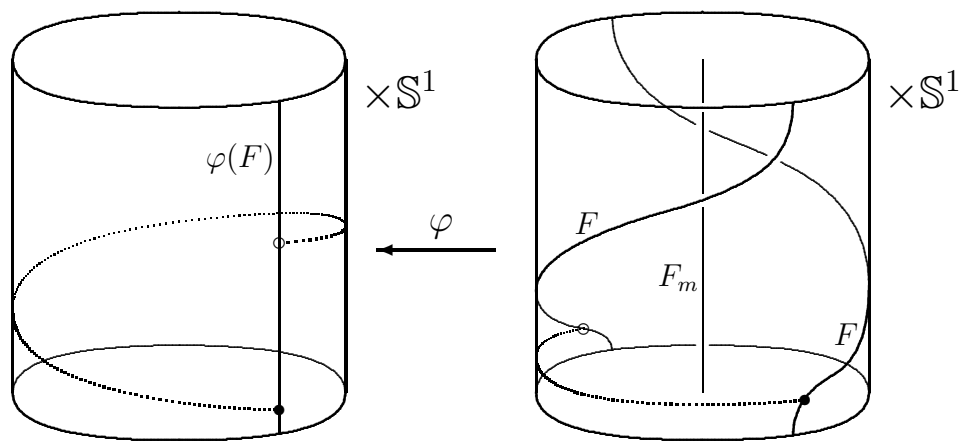
\begin{figure}[bth]
\setlength{\unitlength}{0.9pt}
\begin{picture}(366,200)(-183,-100)

\thicklines
\qbezier(180, 70)(180, 90)(110, 90) 
\qbezier(180, 70)(180, 50)(110, 50)
\qbezier( 40, 70)( 40, 90)(110, 90)
\qbezier( 40, 70)( 40, 50)(110, 50)

\thicklines
\qbezier(180,-70)(180,-90)(110,-90) 
\qbezier( 40,-70)( 40,-90)(110,-90)
\thinlines
\qbezier(180,-70)(180,-50)(113,-50)
\qbezier( 40,-70)( 40,-50)(107,-50)

\thicklines
\put(180,70){\line(0,-1){140}} 
\put( 40,70){\line(0,-1){140}}

\thinlines
\put(110,70){\line(0,-1){18}} 
\put(110,48){\line(0,-1){34}}
\put(110,8){\line(0,-1){78}}

\thicklines
\qbezier(148,52)(148,32)(133,22) 
\qbezier(133,22)(118,12)(90,5)
\qbezier(90,5)(40,-10)(40,-27)
\thinlines
\qbezier(40,-27)(40,-40)(60,-43)
\qbezier(60,-43)(72,-45)(72,-52)

\thicklines
\qbezier(148,-88)(148,-72)(164,-63) 
\qbezier(164,-63)(180,-55)(180,-25)
\thinlines
\qbezier(180,-25)(180,10)(143,27)
\qbezier(137,30)(120,37)(113,40)
\qbezier(107,43)(87,53)(79,66)
\qbezier(79,66)(72,77)(72,88)

\thinlines
\put(153,-73){\circle*{5}}
\put(60,-43){\circle{5}}

\thicklines
\qbezier[20](153,-73)(132,-75)(110,-75) 
\qbezier[50](110,-75)(40,-75)(40,-55)
\thinlines
\qbezier[25](40,-55)(40,-46)(60,-43)

\put(55,-3){$F$}
\put(164,-50){$F$}
\put(90,-25){$F_{m}$}

\thicklines
\qbezier(-180, 70)(-180, 90)(-110, 90) 
\qbezier(-180, 70)(-180, 50)(-110, 50)
\qbezier( -40, 70)( -40, 90)(-110, 90)
\qbezier( -40, 70)( -40, 50)(-110, 50)

\thicklines
\qbezier( -180,-70)( -180,-90)( -110,-90) 
\qbezier(  -40,-70)(  -40,-90)( -110,-90)
\thinlines
\qbezier( -180,-70)( -180,-50)( -110,-50)
\qbezier(  -40,-70)(  -40,-50)( -110,-50)

\thicklines
\put(-180,70){\line(0,-1){140}}
\put( -40,70){\line(0,-1){140}}

\thicklines
\put(-68,53){\line(0,-1){140}} 

\thinlines
\put(-68,-77){\circle*{5}}
\put(-68,-7){\circle{5}}

\thicklines
\qbezier[20](-68,-77)(-89,-77)(-110,-74) 
\qbezier[50](-110,-74)(-180,-64)(-180,-35)
\thinlines
\qbezier[50](-180,-35)(-180,-7)(-110,2)
\qbezier[50](-110,2)(-40,10)(-40,0)
\thicklines
\qbezier[15](-40,0)(-40,-5)(-68,-7)

\put(-100,25){$\phi(F)$} 

\thicklines
\put(23,-10){\vector(-1,0){50}} 
\put(0,0){\makebox(0,0){\large$\phi$}}
\put(185,50){{\Large $\times\Sph{1}$}}
\put(-35,50){{\Large $\times\Sph{1}$}}

\end{picture}
\caption{Logarithmic transformation (multiplicity 2)}
\label{fig-log}
\end{figure}
}

\begin{piece}
Another surgery is the \emph{logarithmic transformation} performed along a torus $T$ that is transverse to a foliation and has $T\cdot T=0$. One takes a neighborhood $T\times D^2$ of the torus $T$, cuts it out and glues it back in using an automorphism of the boundary $\del(T\times D^2)=\Sph{1}\times\Sph{1}\times\Sph{1}$ (see \Figref{fig-log}, where a logarithmic transformation of multiplicity $2$ is suggested, with the unrepresented $\Sph{1}$-factor belonging to $T$). 
If $T$ is transverse to the foliation $\F$, then one way of continuing the foliation across the glued-in $T\times D^2$ is suggested in \Figref{fig-log.filled}, but it needs a (circle of) Reeb components and thus ruins tautness.
\end{piece}

\comment{
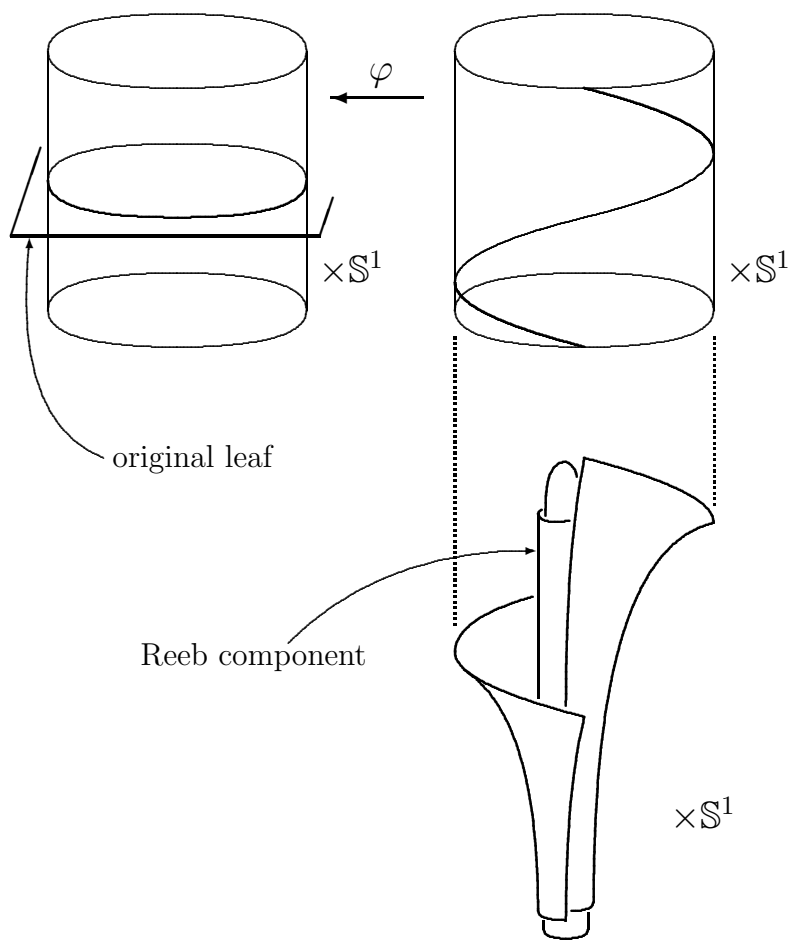
\begin{figure}[bth]
\setlength{\unitlength}{0.7pt}
\begin{picture}(366,500)(-183,-400)

\thinlines
\qbezier(180, 70)(180, 90)(110, 90) 
\qbezier(180, 70)(180, 50)(110, 50)
\qbezier( 40, 70)( 40, 90)(110, 90)
\qbezier( 40, 70)( 40, 50)(110, 50)

\thinlines
\qbezier(180,-70)(180,-90)(110,-90) 
\qbezier( 40,-70)( 40,-90)(110,-90)
\thinlines
\qbezier(180,-70)(180,-50)(110,-50)
\qbezier( 40,-70)( 40,-50)(110,-50)

\thinlines
\put(180,70){\line(0,-1){140}} 
\put( 40,70){\line(0,-1){140}}

\thicklines
\qbezier(180,15)(180,32)(110,50) 
\qbezier( 40,-55)( 40,-72)(110,-90)
\thinlines
\qbezier(180,15)(180,-3)(110,-20)
\qbezier( 40,-55)( 40,-37)(110,-20)

\thinlines
\qbezier(-180, 70)(-180, 90)(-110, 90) 
\qbezier(-180, 70)(-180, 50)(-110, 50)
\qbezier( -40, 70)( -40, 90)(-110, 90)
\qbezier( -40, 70)( -40, 50)(-110, 50)

\thinlines
\qbezier( -180,-70)( -180,-90)( -110,-90) 
\qbezier(  -40,-70)(  -40,-90)( -110,-90)
\thinlines
\qbezier( -180,-70)( -180,-50)( -110,-50)
\qbezier(  -40,-70)(  -40,-50)( -110,-50)

\thinlines
\put(-180,70){\line(0,-1){140}}
\put( -40,70){\line(0,-1){140}}

\thicklines
\qbezier( -180,0)( -180,-20)( -110,-20) 
\qbezier(  -40,0)(  -40,-20)( -110,-20)
\thinlines
\qbezier( -180,0)( -180,20)( -110,20)
\qbezier(  -40,0)(  -40,20)( -110,20)

\thicklines
\put(-200,-30){\line(1,0){167}}
\put(-200,-30){\line(1,3){16}}
\put(-33,-30){\line(1,3){7}}

\thicklines
\qbezier(180,15)(180,32)(110,50) 
\qbezier( 40,-55)( 40,-72)(110,-90)
\thinlines
\qbezier(180,15)(180,-3)(110,-20)
\qbezier( 40,-55)( 40,-37)(110,-20)

\thicklines
\qbezier( 40,-255)( 40,-272)(110,-290) 

\qbezier( 40,-255)( 40,-237)(82,-225)

\qbezier(180,-185)(180,-168)(110,-150) 

\qbezier(110,-150)(100,-200)(100,-284)
\qbezier(180,-185)(115,-200)(115,-390)
\qbezier(115,-390)(115,-395)(103,-395)

\qbezier(110,-290)(100,-330)(100,-400)
\qbezier(45,-264.5)(85,-290)(85,-395)
\qbezier(85,-395)(85,-400)(100,-400)

\put(85,-280){\line(0,1){100}} 
\qbezier(85,-180)(85,-186)(102,-184)
\qbezier(85,-180)(85,-179)(87,-178)

\qbezier(89,-181)(88,-152)(100,-152) 
\qbezier(100,-152)(105,-154)(106,-160)

\put(88,-400){\line(0,-1){5}} 
\put(112,-396){\line(0,-1){9}}
\qbezier(88,-405)(88,-410)(100,-410)
\qbezier(112,-405)(112,-410)(100,-410)

\thicklines
\qbezier[30](180,-175)(180,-130)(180,-85)
\qbezier[52](40,-240)(40,-162)(40,-85)

\thinlines
\put(84,-200){\vector(1,0){0}}
\qbezier(82,-200)(0,-200)(-50,-250)
\put(-130,-262){Reeb component}

\put(-190,-31){\vector(0,1){0}}
\qbezier(-190,-33)(-200,-130)(-150,-150)
\put(-145,-155){original leaf}

\thicklines
\put(23,45){\vector(-1,0){50}} 
\put(0,57){\makebox(0,0){\large$\phi$}}
\put(185,-55){{\large $\times\Sph{1}$}}
\put(-35,-55){{\large $\times\Sph{1}$}}
\put(155,-350){{\large$\times\Sph{1}$}}

\end{picture}
\caption{Filling-in a logarithmic transformation}
\label{fig-log.filled}
\end{figure}
}


\providecommand{\bysame}{\leavevmode\hbox to3em{\hrulefill}\thinspace}
\providecommand{\MR}{\relax\ifhmode\unskip\space\fi MR }
\providecommand{\MRhref}[2]{%
  \href{http://www.ams.org/mathscinet-getitem?mr=#1}{#2}
}
\providecommand{\href}[2]{#2}


\begin{thebibliography}{Rum79}

\bibitem[AL94]{audin.lafontaine}
Mich{\`e}le Audin and Jacques Lafontaine (eds.), \emph{Holomorphic curves in
  symplectic geometry}, Progress in Mathematics, vol. 117, Birkh\"auser Verlag,
  Basel, 1994. \MR{95i:58005}

\bibitem[Bri84]{fabiano}
Fabiano Brito, \emph{A remark on minimal foliations of codimension two}, Tohoku
  Math. J. (2) \textbf{36} (1984), no.~3, 341--350. \MR{85m:53029}

\bibitem[CC00]{candelconlon}
Alberto Candel and Lawrence Conlon, \emph{Foliations. {I}}, Graduate Studies in
  Mathematics, vol.~23, American Mathematical Society, Providence, RI, 2000.
  \MR{2002f:57058}

\bibitem[Don99]{donaldson.lefschetz}
S.~K. Donaldson, \emph{Lefschetz pencils on symplectic manifolds}, J.
  Differential Geom. \textbf{53} (1999), no.~2, 205--236. \MR{2002g:53154}

\bibitem[EM98]{eliashberg}
Y.~Eliashberg and N.~M. Mishachev, \emph{Wrinkling of smooth mappings.
  {I}{I}{I}. {F}oliations of codimension greater than one}, Topol. Methods
  Nonlinear Anal. \textbf{11} (1998), no.~2, 321--350. \MR{2000a:57073}

\bibitem[Gab83]{gabai}
David Gabai, \emph{Foliations and the topology of $3$-manifolds}, J.
  Differential Geom. \textbf{18} (1983), no.~3, 445--503. \MR{86a:57009}

\bibitem[Gab01]{gabai.survey}
\bysame, \emph{3 lectures on foliations and laminations on 3-manifolds},
  Laminations and foliations in dynamics, geometry and topology (Stony Brook,
  NY, 1998), Contemp. Math., vol. 269, Amer. Math. Soc., Providence, RI, 2001,
  pp.~87--109. \MR{2002g:57032}

\bibitem[GN95]{garcianaveira}
F.~J. Garc{\'{\i}}a and A.~M. Naveira, \emph{Two remarks about foliations and
  minimal foliations of codimension greater than two}, Analysis and geometry in
  foliated manifolds (Santiago de Compostela, 1994), World Sci. Publishing,
  River Edge, NJ, 1995, pp.~29--38. \MR{97g:53034}

\bibitem[Gro85]{gromov.pseudo}
M.~Gromov, \emph{Pseudoholomorphic curves in symplectic manifolds}, Invent.
  Math. \textbf{82} (1985), no.~2, 307--347. \MR{87j:53053}

\bibitem[GS99]{gompf}
Robert~E. Gompf and Andr{\'a}s~I. Stipsicz, \emph{$4$-manifolds and {K}irby
  calculus}, American Mathematical Society, Providence, RI, 1999.
  \MR{2000h:57038}

\bibitem[HL82]{harveylawson}
Reese Harvey and H.~Blaine Lawson, Jr., \emph{Calibrated foliations (foliations
  and mass-minimizing currents)}, Amer. J. Math. \textbf{104} (1982), no.~3,
  607--633. \MR{84h:53095}

\bibitem[KM94]{kronheimer.mrowka.proj}
P.~B. Kronheimer and T.~S. Mrowka, \emph{The genus of embedded surfaces in the
  projective plane}, Math. Res. Lett. \textbf{1} (1994), no.~6, 797--808.
  \MR{96a:57073}

\bibitem[Kro99]{kronheimer}
P.~B. Kronheimer, \emph{Minimal genus in ${S}\sp 1\times {M}\sp 3$}, Invent.
  Math. \textbf{135} (1999), no.~1, 45--61. \MR{2000c:57071}

\bibitem[OS00]{ozsvath.szabo.sympl.thom}
Peter Ozsv{\'a}th and Zolt{\'a}n Szab{\'o}, \emph{The symplectic {T}hom
  conjecture}, Ann. of Math. (2) \textbf{151} (2000), no.~1, 93--124.
  \MR{2001a:57049}

\bibitem[Rum79]{rummler}
Hansklaus Rummler, \emph{Quelques notions simples en g\'eom\'etrie riemannienne
  et leurs applications aux feuilletages compacts}, Comment. Math. Helv.
  \textbf{54} (1979), no.~2, 224--239. \MR{80m:57021}

\bibitem[Sco03]{scorpan.foliation.exist}
Alexandru Scorpan, \emph{Existence of foliations on 4-manifolds}, arXiv:
  math.GT/0302318.

\bibitem[Sul76]{sullivan-cycles}
Dennis Sullivan, \emph{Cycles for the dynamical study of foliated manifolds and
  complex manifolds}, Invent. Math. \textbf{36} (1976), 225--255. \MR{55
  \#6440}

\bibitem[Sul79]{sullivan-taut}
\bysame, \emph{A homological characterization of foliations consisting of
  minimal surfaces}, Comment. Math. Helv. \textbf{54} (1979), no.~2, 218--223.
  \MR{80m:57022}

\bibitem[Thu74]{thurston}
William Thurston, \emph{The theory of foliations of codimension greater than
  one}, Comment. Math. Helv. \textbf{49} (1974), 214--231. \MR{51 \#6846}

\bibitem[Thu76]{thurston-codim1}
W.~P. Thurston, \emph{Existence of codimension-one foliations}, Ann. of Math.
  (2) \textbf{104} (1976), no.~2, 249--268. \MR{54 \#13934}

\bibitem[Thu86]{thurston-norm}
William~P. Thurston, \emph{A norm for the homology of {$3$}-manifolds}, Mem.
  Amer. Math. Soc. \textbf{59} (1986), no.~339, i--vi and 99--130.
  \MR{88h:57014}

\end{thebibliography}

\end{document}